\documentclass[11pt]{article}

\usepackage[english]{babel}

\usepackage[letterpaper,top=2cm,bottom=2cm,left=3cm,right=3cm,marginparwidth=1.75cm]{geometry}

\usepackage{amsmath}
\usepackage{graphicx}
\usepackage[colorlinks=true, allcolors=blue]{hyperref}
\usepackage{multirow}
\usepackage{array}
\usepackage{epsfig,rotating}
\usepackage{amssymb,amsmath, mathabx}
\usepackage{latexsym}
\usepackage[usenames]{color}
\usepackage{upgreek} 

\numberwithin{equation}{section}

\newtheorem{theorem}{Theorem}
\newtheorem{lemma}{Lemma}[section]

\newtheorem{prop}{Proposition}

\newcommand{\be}{\begin{equation}}
\newcommand{\ee}{\end{equation}}
\newcommand{\beaa}{\begin{eqnarray*}}
\newcommand{\eeaa}{\end{eqnarray*}}
\newcommand{\bea}{\begin{eqnarray}}
\newcommand{\eea}{\end{eqnarray}}

\newcommand{\bei}{\begin{itemize}}
\newcommand{\eei}{\end{itemize}}

\def\E{\mathrm{E}}
\def\Var{\mathrm{Var}}

\def\eop{\hfill $\Box$\\}
\def\bed{\begin{description}}\def\eed{\end{description}}
\def\ben{\begin{enumerate}}\def\een{\end{enumerate}}
\def\bea{\begin{eqnarray}}\def\eea{\end{eqnarray}}
\def\bean{\begin{eqnarray*}}\def\eean{\end{eqnarray*}}
\def\ba{\begin{array}}\def\ea{\end{array}}
\def\bt{\begin{theorem}  {\bf }}\def\et{\end{theorem}}
\def\bl{\begin{lemma}{\bf }}\def\el{\end{lemma}}
\def\br{\begin{remark}$\!\!\!\!\!$ {\bf .}}\def\er{\end{remark}}
\def\bc{\begin{corollary}$\!\!\!\!\!$ {\bf }}\def\ec{\end{corollary}}

\def\bz{\mathbf{z}}

\def\P{\mathbb{P}}
\def\E{\mathbb{E}}
\def\Var{\mathrm{Var}}

\begin{document}
    
\title{\sc Eigenvalues of Product of
Ginibre Ensembles and Their Inverses and that of Truncated Haar Unitary Matrices and  Their
Inverses} 



\author{\sc Shuhua Chang, Tiefeng Jiang and Yongcheng Qi}





\date{}


\maketitle

\begin{abstract}
Consider two types of products of independent random matrices,
including products of Ginibre matrices and inverse Ginibre matrices
and products of truncated Haar unitary matrices and inverse truncated
Haar matrices. Each product matrix has $m$ multiplicands of
$n$ by $n$ square matrices, and the empirical distribution based on
the $n$ eigenvalues of the product matrix is called empirical
spectral distribution of the matrix. In this paper, we investigate
the limiting empirical spectral distribution of the product matrices
when $n$ tends to infinity and $m$ changes with $n$. For properly
scaled eigenvalues for two types of the product matrices, we
obtain the necessary and sufficient conditions for the convergence
of the empirical spectral distributions.
\end{abstract}

\section{Introduction}

The study of non-Hermitian matrices was initiated by
Ginibre~\cite{Ginibre1965} who investigated the limiting empirical
spectral distributions of ensembles of real, complex and quaternion matrices with Gaussian probability distributions. The
circular law problems for non-Hermitian matrices had been studied for many
years and were established  by the cornerstone work of
Girko~\cite{Girko, Girko2004}, Bai~\cite{Bai}, Bai and
Silverstein~\cite{BS2006}, and Tao and Vu~\cite{Tao2008, Tao}.  For
more work, see also Pan and Zhou~\cite{Pan}, G\"{o}tze and
Tikhomirov~\cite{GFTA} and literature therein. In this paper, we will focus on the spectral distributions of products of
random matrices.

In recent years products of random matrices have been studied actively. Their applications include  wireless telecommunication, disordered spin chain, the stability of large
complex system, quantum transport in disordered wires, symplectic
maps and Hamiltonian mechanics, quantum chromo-dynamics at non-zero
chemical potentials; see, e.g., Ipsen~\cite{Ipsen2015}. More
information on early study can also be found in Crisanti {\it et al.}~\cite{CPV1993}, Akemann and Ipsen~\cite{AI2015} and
references therein.

For product of $m$ independent matrices, early papers focus on the empirical spectral distributions for the product matrices when $m$ is a fixed integer.
For some universality results under moment
conditions for entries of random matrices; see, e.g., G\"{o}tze and
Tikhomirov~\cite{Goetz}, Bordenave~\cite{Bor}, O'Rourke and
Soshnikov~\cite{Rourke}, O'Rourke {\it et al.}~\cite{Rourke14}. For
various results concerning products of random matrices with Gaussian
entries, see, e.g., Burda {\it et al.}~\cite{BJW}, Burda~\cite{Burda},
Forrester~\cite{Forrester2013, Forrester2014},
Ipsen~\cite{Ipsen2015}, Forrester and Liu~\cite{FL2016}, Liu and
Wang~\cite{LW2016}, Liu {\it et al.}~\cite{LWZ2016}, and Akemann {\it et al.}~\cite{A1, A2, A3, A4, A5, A6, ABKN}. As the
eigenvalues of a random matrix form a determinantal point process, a very recent
paper by Jiang and Qi~\cite{JiangQi2019} investigated
the limiting empirical spectral distribution of the product of $m$
independent $n\times n$ Ginibre matrices and truncated Haar unitary
matrices when $m$ is an arbitrary sequence of positive integers. Jiang and Qi~\cite{JiangQi2017} also obtained the limiting
spectral radii of product matrices for Ginibre ensembles when $m/n$ has a limit in $[0, \infty]$.
Subsequently, for the products of $m$ independent $n\times n$ spherical
ensembles, Zeng~\cite{Zeng2016}
studied the empirical spectral distribution of the product when $m$ is fixed, Chang and Qi~\cite{ChangQi2017}
obtained the limiting empirical spectral distribution when $m$ is an arbitrary sequence of positive integers,
 and Chang {\it et al.} \cite{ChangLiQi2018} investigated the limiting spectral
radius when $m$ is fixed or $m$ tends to infinity as $n$ goes to infinity.
Qi and Zhao~\cite{QiZhao2021} and Qi and
Xie~\cite{QiXie2020} obtained the limiting empirical spectral
distribution and spectral radii for the product of $m$ independent
rectangular Ginibre matrices by allowing $m$ to change with the
dimension of the product matrices in an arbitrary way, and Zeng~\cite{Zeng2017}
considered the limiting empirical distribution when $m$ is fixed.
Wang~\cite{Wang2018} obtained the limiting distributions for the
largest moduli of the eigenvalues of product random matrices from $m$
polynomial ensembles when $m/n$ has a finite limit.

In this paper, we consider  two types of product
matrices. The first one is the product of Ginibre matrices and inverse Ginibre
matrices. The second is the product of truncated Haar unitary matrices and inverse
truncated Haar matrices. More precisely,  let $A_1, \cdots, A_m$ be
$m$ independent $n\times n$ Ginibre matrices, or  $n\times n$
truncated Haar unitary matrices, and define the product matrix
$A_1^{\varepsilon_1}A_2^{\varepsilon_2}\cdots A_m^{\varepsilon_m}$,
where $\{\varepsilon_j, 1\le j\le m\}$ are constants taking values
$1$ or $-1$.
{    These types of product matrices include random matrices which have been well studied in the literature. For example, when $A_j$'s are independent Ginibre matrices,  $A^{-1}_1
A_2$ is the spherical ensemble, and a product of $k$ independent matrices from the spherical ensemble can be written as $A_1^{\varepsilon_1}A_2^{\varepsilon_2}\cdots A_m^{\varepsilon_m}$, where $\varepsilon_j=(-1)^j$ for $1\le j\le m$ with $m=2k$. See more details in Remark 2.} In this  paper, we will investigate
the limiting empirical spectral
distributions of the product matrices when $m$ changes with $n$.  Our empirical spectral distributions are based on
the arguments and scaled moduli of the eigenvalues of product
matrices.

{      Now we present some features of the paper. First, one of our results establishes the limiting law of the empirical distributions of the eigenvalues of a generalized spherical ensemble, a notion defined below \eqref{JQYC}. Furthermore, some connections between our results with  Gaussian free fields, Gaussian multiplicative chaos and quantum computing are established, see Remarks 6 and 7 from Section \ref{PtHum}. Third, we discover a counterintuitive phenomenon appeared in our main results: the limiting empirical distributions of the eigenvalues of
$A_1^{\varepsilon_1}A_2^{\varepsilon_2}\cdots A_m^{\varepsilon_m}$ do not depend the locations of ``$+1$" or ``$-1$", but the limits depend on how many ``$+1$" and ``$-1$" in total, respectively. It is a bit bizarre since the product of matrices are not commentative. See Theorems \ref{ginibre} and \ref{harr}.

Finally, we make some comments about the technical part. We  utilize  some tools developed by Jiang and Qi~\cite{JiangQi2019} for determinantal point processes.
First,  we will give distributional representations for eigenvalues of the products for both Ginibre matrices and truncated Haar unitary matrices and use them to approach the limiting spectral distributions.  In fact, the limiting spectral distributions in this paper are much more complicated than those for the product matrices studied in the literature when inverse matrices are not involved in the products. Second, we give necessary and sufficient conditions for convergence of the empirical spectral distributions in order to  give a complete description of the limiting spectral distributions for both the product matrices. Our approach is more challenging than just giving a sufficient condition for convergence of the empirical spectral distributions. To show the necessities,
we construct some analytic function and then  apply subsequence arguments and use properties of analytic functions.
}

The rest of the paper is organized as follows. In
Section~\ref{main}, we present the main results on limiting
empirical spectral distributions of the two types of product matrices. In
Section~\ref{proofs}, we introduce some auxiliary lemmas and prove
the main results.

\section{Product matrices and main results}\label{main}

In this paper, we assume $A_1, \cdots, A_m$ are independent
$n\times n$ random matrices and $\{\varepsilon_j, 1\le j\le m\}$ are
constants taking values $1$ or $-1$. The product matrix of interest is defined by
\begin{equation}\label{Zm}
Z_m:=A_1^{\varepsilon_1}A_2^{\varepsilon_2}\cdots
A_m^{\varepsilon_m},
\end{equation}
where $A_i^{-1}$ stands for the inverse of matrix $A_i$.
Let $\{\mathbf{z}_j;\, 1\leq j \leq n\}$ be the eigenvalues of $Z_m$. Note that $\{\mathbf{z}_j;\, 1\leq j \leq n\}$
are complex random variables.
Write
 \begin{equation}\label{argument}
\Theta_j=\arg(\mathbf{z}_j)\in [0, 2\pi)~\mbox{ such that }
~\mathbf{z}_j=|\mathbf{z}_j|\cdot e^{i\Theta_j},
\end{equation}
where $i=\sqrt{-1}$ is the imaginary unit. For real-valued function $h_n(x)$
defined on $[0, \infty)$, define
\begin{equation}\label{mun}
\mu_n=\frac{1}{n}\sum^n_{j=1}\delta_{(\Theta_j,
h_n(|\mathbf{z}_j|))}
\end{equation}
as the empirical measure of the pairs $(\Theta_1, h_n(|\mathbf{z}_1|), \cdots, (\Theta_n,
h_n(|\mathbf{z}_n|)$. Similarly, for non-negative random variables $Y_1, \cdots, Y_n$, we
define the empirical measure
\begin{equation}\label{nun}
\nu_n=\frac{1}{n}\sum^n_{j=1}\delta_{h_n(Y_j)}.
\end{equation}
In the paper,  $\{Y_j, ~1\le j\le n\}$ are independent random variables such that $g(Y_1, \cdots, Y_n)$ and $g(|\bz_1|, \cdots, |\bz_n|)$ are identically distributed for any $n$-variate symmetric function $g$. Certainly, $\mu_n$ depends on the pairs of $(\Theta_j, h_n(|\mathbf{z}_j|), 1\leq j \leq n$ and $\nu_n$ depends on random variables $Y_1, \cdots, Y_n$, respectively. There will be no potential confusion from the context in our discussions later on.

A special attention is paid to the case that $h_n(x)$ is a linear function. In this situation
the corresponding empirical spectral measure is especially denoted by
\begin{equation}\label{mu*}
\mu_n^*=\frac1n\sum^n_{j=1}\delta_{\mathbf{z}_j/a_n}\ \ \  \mbox{and}\ \  \ \nu_n^*=\frac{1}{n}\sum^n_{j=1}\delta_{Y_j/a_n}
\end{equation}
where $a_n>0$ is a sequence of normalizing constants. The choice of $a_n$ is case by case.

The study of convergence of $\mu_n$ given in \eqref{mun} is
initiated by Jiang and Qi~\cite{JiangQi2019}. As $m$ grows with
$n$,  the moduli for eigenvalues of a product matrix usually grow to
infinity or vanish to zero exponentially fast,  and the measure
$\mu_n^*$ may not converge to a measure for any sequence $\{a_n\}$.
By choosing a proper nonlinear transformation $h_n(\cdot)$ for the moduli, one
can obtain a non-singular limiting measure of $\mu_n$. These are seen in the definitions of \eqref{mun} and \eqref{nun}.

With the notation in \eqref{argument}, we may use $(\Theta_j,
h_n(\bz_j))$ to form a new complex number
$h_n(|\bz_j|)e^{i\Theta_j}$.  Therefore, we can define the empirical
spectral measure for scaled eigenvalues $h_n(|\bz_j|)e^{i\Theta_j}$
such that
\begin{equation}\label{hatmu}
\hat\mu_n=\frac1n\sum^n_{j=1}\delta_{h_n(|\mathbf{z}_j|)e^{i\Theta_j}}.
\end{equation}
It is easy to see that the two measures $\hat\mu_n$ and $\mu_n^*$ are
the same when $h_n(r)=r/a_n$. For clarity, in this case we use $\mu_n^*$ to denote the empirical measure of
$\mathbf{z}_j$'s in the future.

For a sequence of random probability measures $\{\tau, \tau_n;\,
n\geq 1\}$, we write
 \[
\tau_n \rightsquigarrow \tau\ \ \mbox{if \ $\mathbb{P}$($\tau_n$
converges weakly to $\tau$ as $n\to\infty$)=1}. \]
 When $\tau$ is a
non-random probability measure generated by random variable $X$, we
sometimes simply write $\tau_n \rightsquigarrow X$.
Throughout the paper, Unif($A$) denotes the uniform distribution on
a set $A$.

 Recall $A_1, \cdots, A_m$ and $Z_m$ defined in \eqref{Zm}. In this paper, we will consider two different types of matrices $A_i$'s. They will be specified later and the corresponding limiting spectral distribution of $Z_m$ for each case will be investigated.  A
common feature for the two product matrices is that the eigenvalues
of both form a determinantal point process on complex plane
$\mathbb{C}$.

\subsection{Product of
Ginibre matrices and their inverses}

A standard complex normal distribution $\mathbb{C}N(0, 1)$ is the distribution of $(\xi_1+i\xi_2)/\sqrt{2}$ where $\xi_1$ and $\xi_2$ are i.i.d. $N(0, 1)$-distributed random variables.

Assume $A_1, \cdots, A_m$ are independent $n\times n$ random
matrices whose $mn^2$ entries are i.i.d. complex normal random
variables with distribution $\mathbb{C}N(0, 1)$. They are called complex Ginibre matrices.
Assume also $\{\varepsilon_j, 1\le j\le m\}$ are constants taking
values $1$ or $-1$.
Each Ginibre matrix is invertible since its eigenvalues have a continuous joint density function;  see, e.g., Ginibre~\cite{Ginibre1965}.


For the joint density function of the eigenvalues $\bz_1, \cdots,
\bz_n$ of $Z_m$, we cite a known result which relates $\bz_1,
\cdots, \bz_n$ to a rotation-invariant determinantal point process and plays a key role in obtaining the limiting expected empirical spectral distributions
in Adhikari {\it et al.}~\cite{ARRS2016} when $m$ is a fixed integer.

\begin{prop}\label{prop1} (Theorem 1 from Adhikari {\it et al.}~\cite{ARRS2016}) Suppose $A_1, \cdots, A_m$ are independent and identically distributed (i.i.d.) $n\times n$ random
matrices and the $n^2$ entries of $A_1$ are i.i.d. random variables of distribution $\mathbb{C}N(0, 1)$. Let
$\{\varepsilon_j, 1\le j\le m\}$ be constants taking values $1$ or
$-1$. Let $\bz_1, \cdots, \bz_n$ be the eigenvalues of
$Z_m=A_1^{\varepsilon_1}A_2^{\varepsilon_2}\cdots
A_m^{\varepsilon_m}$. Then,  $\bz_1, \cdots, \bz_n$ form a
determinantal point process with kernel
\[
\mathbb{K}_n(z,w)=\sqrt{\varphi(z)\varphi(w)}\sum^{n-1}_{r=0}\frac{(z\bar{w})^r}
 {(2\pi)^m(r!)^p\big[(n-r-1)!\big]^{m-p}}
\]
with respect to the Lebesgue measure on $\mathbb{C}$, where
\begin{equation}\label{p}
    p=\#\{1\leq i \leq m: \varepsilon_i=1\},
\end{equation}
and $\varphi(z)$ is a
weight function such that
\begin{equation}\label{dphi}
\varphi(z)dz=\int_{z_1^{\varepsilon_1}\cdots
z_m^{\varepsilon_m}=z}e^{-\sum^m_{j=1}|z_j|^2}\prod^m_{j=1}|z_j|^{(1-\varepsilon_j)(n-1)}\prod^m_{j=1}dz_j
\end{equation}
and  $dz$ is the Lebesgue measure on complex plane and equivalent to
$dz=dxdy$ with $z=x+yi$ on $\mathbb{R}^2$ space.  Equivalently, $\bz_1,
\cdots, \bz_n$ have a joint probability density function
proportional to
\[
\prod_{1\le j<k\le n} |z_j-z_k|^2\prod^n_{\ell=1}\varphi(z_\ell)
\]
with respect to the Lebesgue measure on $\mathbb{C}^n$.
\end{prop}

An explicit expression for $\varphi$ seems difficult to obtain from
\eqref{dphi}. But the background measure $\varphi(z)dz$ is
rotation-invariant and $\varphi(z)=\varphi(|z|)$; see also Remark 10
in Adhikari {\it et al.}~\cite{ARRS2016}.  Moreover, one should interpret
\[
\varphi(z)dz=\int_{h(z_1,\cdots, z_m)=z}g(z_1,\cdots,
z_n)\prod^m_{j=1}dz_j
\]
as
\begin{equation}\label{rep}
\int f(z)\varphi(z)dz=\int f(h(z_1,\cdots, z_m))g(z_1,\cdots,
z_n)\prod^m_{j=1}dz_j
\end{equation}
for all $f$: $\mathbb{C}\to \mathbb{C}$ integrable functions.   We
will explore $\varphi$ through identity \eqref{rep}. For details,
see Lemma~\ref{t-moment}.

 A probabilistic interpretation of $\varphi$ under \eqref{dphi} is that $\varphi(z)$ is proportional to the density function of the complex random variable $\bz_1^{\varepsilon_1}\cdots
\bz_m^{\varepsilon_m}$ where $\bz_1, \cdots, \bz_m$ have joint density function proportional to
\beaa
e^{-\sum^m_{j=1}|z_j|^2}\prod^m_{j=1}|z_j|^{(1-\varepsilon_j)(n-1)}
\eeaa
with respect to the Lebesgue measure on $\mathbb{R}^{2m}$, where we understand $z_j=x_j+y_ji$ for each $j$ and $(z_1, \cdots, z_m)$ as $(x_1, y_1,  x_2,  y_2, \cdots, x_m, y_m)$. Obviously, $\bz_1, \cdots, \bz_m$ are independent random variables and $z_j$ has a density function proportional to $e^{-|z|^2}|z|^{(1-\epsilon_j)(n-1)}$ for each $j$.
In particular, when $m=1$ and $\varepsilon_1=1$, $\varphi(z)=\exp(-|z|^2)$, and $Z_1$ is the Ginibre ensemble. When $m=2$, $\varepsilon_1=-1$ and $\varepsilon_2=1$, we have $\varphi(z)=n!\pi({1+|z|^2})^{-(n+1)}$, and in this case, $Z_1=A_1^{-1}A_2$ is the spherical ensemble.

For a better understanding on determinantal point processes,
interested readers can see, for example, Soshnikov~\cite{SOS2000},
Johansson~\cite{Jo2005} and Hough {\it et al.}~\cite{HKPV2009}. When $m$ is a fixed integer Adhikari {\it et al.}~\cite{ARRS2016} obtained
the limiting expected empirical distribution for square radii of
eigenvalues of $Z_m$. In this paper we consider the limiting spectral distribution
of the eigenvalues and allow $m$ to change with $n$. To state our results we need some notation.
For any $\alpha\in [0,1]$, define function
\[
G_\alpha(x)=\frac{x^{\alpha}}{(1-x)^{1-\alpha}},~~~x\in (0,1),
\]
and  $G_\alpha(0)=G_\alpha(0+)=\lim_{x\downarrow 0}G_{\alpha}(x)$
and $G_\alpha(1)=G_\alpha(1-)=\lim_{x\uparrow 1}G_{\alpha}(x)$. Then
$G_\alpha(x)$ is continuous and strictly increasing in $[0,1)$. It is easy to verify that
\begin{equation}\label{boundofG}
x\le G_{\alpha}(x)\le \frac{1}{1-x} ~~~\mbox{and}~~~ 0<G_{\alpha}'(x)\le
\frac{1}{x(1-x)^2}, ~~~x\in (0,1)
\end{equation}
 for any $\alpha\in [0,1]$. Again, for each $\alpha\in [0,1]$, define
\begin{equation}\label{G*}
G^*_\alpha(x)=\left\{
         \begin{array}{ll}
           0, & \hbox{ if } x\le G_\alpha(0); \\
           G^{-1}_\alpha(x), & \hbox{ if } x\in (G_\alpha(0), G_\alpha(1));  \\
           1, & \hbox{ if } G_\alpha(1)\le x<\infty.
         \end{array}
       \right.
\end{equation}
The last case in the definition of $G_{\alpha}^*$ is needed only if
$G_\alpha(1)<\infty$.   For every $\alpha\in [0,1]$,
$G^*_{\alpha}(x)$ is well defined in $\mathbb{R}$, and is
non-decreasing in $x$ with the following property
\begin{equation}\label{propertyof G*}
\lim_{x\downarrow 0}G^*_{\alpha}(x)=0~~\mbox{ and }~~
\lim_{x\uparrow\infty}G^*_{\alpha}(x)=1.
\end{equation}
This concludes that $G^*_{\alpha}(x)$ is a probability
distribution function. In fact, $G^*_\alpha$ is continuous in
$\mathbb{R}$ and differentiable on $\big((G_\alpha(0),
G_\alpha(1)\big)$ with density function
$g^*_\alpha(x)=\frac{d}{dx}G^*_{\alpha}(x)I(G_\alpha(0)<x<
G_\alpha(1))$. Write
$G_{\alpha}(x)=(1-x)^{-1}\big(x(1-x)\big)^{\alpha}$. Evidently,  $G_{\alpha}(x)$ is strictly
decreasing in $\alpha\in [0,1]$, that is,
\begin{equation}\label{OrderofG*}
G^*_0(x)\le G^*_{\alpha}(x)\le G^*_1(x),~~~\alpha\in [0,1] ~,  x\in \mathbb{R}.
\end{equation}

Our main result in this section is as follows.

\begin{theorem}\label{ginibre}  Let $m=m_n$ be a sequence of
positive integers. Set $a_n=n^{2p-m}$ and
$h_n(r)=(r^2/a_n)^{1/\gamma_n}$, $r\ge 0$, where $\{\gamma_n\}$ is a
sequence of positive numbers, and $p=p_n$ is a sequence of non-negative integers as define in \eqref{p}.\\
\noindent{(i)}  If $m_n/\gamma_n\to 0$ as $n\to\infty$,
then $\mu_n \rightsquigarrow \mu=\mathrm{Unif[0,2\pi)}\otimes
\delta_1$, where $\delta_1$ is the delta function at $1$,
representing a degenerate
probability measure at $1$;\\
\noindent{(ii)} If for some $\alpha$ and $\beta$ such that
\begin{equation}\label{twolimits}
\lim_{n\to\infty}\frac{p}{m}=\alpha\in [0,1] \mbox{ and }
\lim_{n\to\infty}\frac{m}{\gamma_n}=\beta\in (0,\infty),
\end{equation}
then $\mu_n \rightsquigarrow \mu=\mathrm{Unif[0,2\pi)}\otimes \nu$,
where the probability measure $\nu$ is determined by
$\nu([0,y])=G^*_{\alpha}(y^{1/\beta})$ for all $y>0$.  Conversely,
if $\mu_n \rightsquigarrow \mu$, where $\mu$ is non-singular
probability measure, then \eqref{twolimits} holds for some $\alpha$
and $\beta$.
\end{theorem}

\noindent\textbf{Remark 1.} A probability measure $\mu$ is singular
if there exists a measurable set $B$ with Lebesgue  measure zero
such that $\mu(B)=1$.

\noindent\textbf{Remark 2.} From Theorem~\ref{ginibre}, we can
obtain some interesting results under special conditions as follows.

\begin{enumerate}
  \item  When $\varepsilon_j=1$ for $1\le j\le m_n$,  $Z_{m_n}$ is a
  product of $m_n$ independent Ginibre matrices,  and the empirical
  spectral distribution of $Z_{m_n}$ was obtained by Jiang and
  Qi~\cite{JiangQi2019}.   In this case,  $\alpha_n=p/m_n=1$.  By setting $\gamma_n=m_n$,
  \eqref{twolimits} holds with $\alpha=1$ and $\beta=1$. Since
  $G_1(x)=x$ for $x\in (0,1)$, and $G_1^*(x)=x$, $x\in (0,1)$.   The
  induced probability measure $\nu$ is \textrm{Unif}($[0,1]$).  Immediately,
  we have
\[
\mu_n=\frac1n\sum^n_{j=1}\delta_{(\Theta_j,
\frac{1}{n}\,|\bz_j|^{2/m_n})}\rightsquigarrow \mathrm{Unif}\big([0,
2\pi)\otimes [0,1]\big)
\]
as $n\to\infty$. If $m_n=m$ is a fixed integer,  we have
\[
\mu_n^*:=\frac1n\sum^n_{j=1}\delta_{\bz_n/n^{m/2}}\rightsquigarrow
\mu^*,
\]
where $\mu^*$ is a probability measure on $\mathcal{C}$ with
density function $\frac{1}{m\pi}|z|^{\frac2m-2}$, $|z|\le 1$.  See
Theorem 2 in  Jiang and Qi~\cite{JiangQi2019} and the statement
following the theorem.

  \item   When $\varepsilon_{2j-1}=-1$ and $\varepsilon_{2j}=1$ for $1\le j\le k_n$, and set
  $m_n=2k_n$.  Then $Z_{2k_n}$ is a product of $k_n$ independent
  and identically distributed spherical matrices.  The limiting empirical
  distribution of eigenvalue of $Z_{2k_n}$ was obtained by Chang and
  Qi~\cite{ChangQi2017}.   By setting $\gamma_n=m_n=2k_n$, we have
  $\alpha=\frac12$ and $\beta=1$. Then our Theorem~\ref{ginibre} leads to
  Theorem~2.1 in Chang and Qi~\cite{ChangQi2017}, that is,
\[
\mu_n=\frac1n\sum^n_{j=1}\delta_{(\Theta_j,
|\bz_j|^{2/m_n})}\rightsquigarrow \mu,
\]
where $\mu$ is a probability measure with density
\[
f(\theta, r)=\frac{1}{\pi}\frac{r}{(1+r^2)^2}, ~~\theta\in [0,2\pi),~
r \in (0,\infty).
\]
If $m_n=2k$ for a fixed integer $k\ge 1$,  then
\[
\mu_n^*:=\frac1n\sum^n_{j=1}\delta_{\bz_n}\rightsquigarrow
\mu^*,
\]
where $\mu^*$ has a density function
\begin{eqnarray}\label{JQYC}
\kappa(z):=\frac{1}{k\pi}\frac{|z|^{1/k-2}}{(1+|z|^{1/k})^2}, z\in \mathbb{C}.
\end{eqnarray}

{      We point out that the convergence of $\hat{\mu}_n$ as defined in \eqref{hatmu} can
also be established by a similar argument employed  in Jiang and Qi~\cite{JiangQi2019}. We skip the details for the length of the paper. }

{      Now we consider a generalized spherical ensemble. Taking   $\varepsilon_1=\cdots =\varepsilon_k=-1$ and $\varepsilon_{k+1}=\cdots =\varepsilon_{2k}=1$ in \eqref{Zm}, the matrix
$Z_{2k}=A_1^{\varepsilon_1}A_2^{\varepsilon_2}\cdots
A_{2k}^{\varepsilon_{2k}}$ becomes $S_1^{-1}S_{2}$, where $S_1$ and $S_2$ are two i.i.d. random matrices and each is a product of $k$ i.i.d. Ginibre  matrices. If $k$ is fixed, then $\mu_n^*$ converges to a probability distribution with density function $\kappa(z)$ defined in \eqref{JQYC}. As $k=1$, the limiting density  $\kappa(z)$ is the same as that in Krishnapur~\cite{Krishnapur2009} for the spherical ensemble. }

\end{enumerate}

\subsection{Product of truncated Haar unitary matrices}\label{PtHum}

Let $U_1, \cdots, U_m$ be $m$ independent Haar distributed unitary
matrices of dimension $n_j\times n_j$ for $j=1, \cdots, m$
respectively, where $n\le n_j$ and $A_1, \cdots, A_m$ be $n\times n$
left uppermost blocks of $U_1, \cdots, U_m$, respectively.  Matrices $A_j$'s are called truncated Haar unitary matrices because each $A_j$ is obtained after we remove the last $n_j-n$ rows and $n_j-n$ columns from $U_j$.
Assume
$\{\varepsilon_j, 1\le j\le m\}$ are constants taking values $1$ or
$-1$.  Then define $Z_m$ as in \eqref{Zm}.
All truncated Haar unitary matrices are invertible because their eigenvalues have a joint density function;  see, e.g., \.{Z}yczkowski and Sommers~\cite{Zski}.

When $m$ is a fixed integer Adhikari {\it et al.}~\cite{ARRS2016} obtained
the limiting expected empirical distribution for square radii of
eigenvalues of $Z_m$. We consider the limiting spectral distribution
of the eigenvalues and allow $m$ to change with $n$.  The
eigenvalues form a determinantal point process as proved in Adhikari
{\it et al.}~\cite{ARRS2016}; see Proposition~\ref{prop2} below.   In
Proposition~\ref{prop2}, we need the beta function $\mathrm{B}(a,b)$ which is defined as
\begin{equation}\label{beta}
\mathrm{B}(a,b)=\int^1_0t^{a-1}(1-t)^{b-1}dt
\end{equation}
for any $a>0$ and $b>0$.

\begin{prop}\label{prop2} (Theorem 3 from Adhikari {\it et al.}~\cite{ARRS2016})  Assume $A_1, \cdots, A_m$ are $m$ truncated
Haar unitary matrices defined above. Then the eigenvalues $\bz_1,
\cdots, \bz_n$ of $Z_m=A_1^{\varepsilon_1}A_2^{\varepsilon_2}\cdots
A_m^{\varepsilon_m}$ form a determinantal point process with kernel
\[
 \mathbb{K}_n(z,w)=\sqrt{\varphi(z)\varphi(w)}\sum^{n-1}_{r=0}\frac{(z\bar{w})^r}{(2\pi)^mC_r}
\]
with respect to the Lebesgue measure on $\mathbb{C}$, where
 \[
 C_r=\prod^m_{k=1}\mathrm{B}\Big(n\frac{1-\varepsilon_k}{2}+\frac{1+\varepsilon_k}{2}+r\varepsilon_k, n_k-n\Big)
 \]
 and
$\varphi(z)$ is a weight function with
\begin{equation}\label{dphi2}
\varphi(z)dz=\int_{z_1^{\varepsilon_1}\cdots
z_m^{\varepsilon_m}=z}\prod^m_{j=1}(1-|z_j|^2)^{n_j-n-1}|z_j|^{(n-1)
(1-\varepsilon_j)}I(|z_j|\le 1)dz_j.
\end{equation}
Equivalently, $\bz_1, \cdots, \bz_n$ have a joint probability
density function proportional to
\[
\prod_{1\le j<k\le n} |z_j-z_k|^2\prod^n_{\ell=1}\varphi(z_\ell)
\]
with respect to the Lebesgue measure on $\mathbb{C}^n$.
\end{prop}

Our main result on limiting empirical spectral distribution  for the
product matrix $Z_m$ is as follows.

\begin{theorem}\label{harr} Define
$a_n=\prod^m_{k=1}(\frac{n}{2n_k-n})^{\varepsilon_k}$,
\begin{equation}\label{deltanj}
\delta_{n,j}=\frac1j\sum^m_{k=1}(-\varepsilon_k)^{j-1}\Big(1-\big(\frac{n}{2n_k-n}\big)^j
\Big),~~~j\ge 1.
\end{equation}
Set $\Delta_n=\delta_{n,1}=\sum^m_{k=1}\frac{2(n_k-n)}{2(n_k-n)+n}$,
and $h_n(r)=(r^2/a_n)^{1/\gamma_n}$, $r\ge 0$, where $\{\gamma_n\}$
is a sequence of positive numbers.\\
\noindent{(i)}
 If $\Delta_n/\gamma_n\to 0$ as $n\to\infty$,
then $\mu_n \rightsquigarrow \mu=\mathrm{Unif[0,2\pi)}\otimes
\delta_1$.

\noindent{(ii)} Assume
$\lim_{n\to\infty}n\Delta_n=\infty$. If
\begin{equation}\label{limits}
\lim_{n\to\infty}\frac{\delta_{n,j}}{\gamma_n}=:\beta_j~\mbox{exists
for each $j\ge 1$, and } \beta_1\in (0,\infty),
\end{equation}
then $\mu_n \rightsquigarrow \mu=\mathrm{Unif[0,2\pi)}\otimes \nu$,
with $\nu$ being a probability measure induced by
$\nu((-\infty,y])=F^*(\log y)$ for all $y>0$, where the probability
cumulative distribution $F^*$ is defined as
 \begin{equation}\label{F*}
F^*(y)=\left\{
         \begin{array}{ll}
           0, & \hbox{ if } y\le f(0), \\
           f^{-1}(y), & \hbox{ if } f(0)<y<f(1), \\
           1, & \hbox{ if } y\ge f(1),
         \end{array}
       \right.
\end{equation}
and $f^{-1}$ is the inverse function of $f$ given by
\begin{equation}\label{ff}
f(x)=\sum^{\infty}_{j=1}\beta_j\big(2x-1\big)^j,~~0<x<1.
\end{equation}
The function $f(x)$ is analytic and strictly increasing on $(0,1)$.

Conversely, if $\mu_n \rightsquigarrow \mu$, where $\mu$ is a
non-singular probability measure, then \eqref{limits} holds.
\end{theorem}

\noindent\textbf{Remark 3.}  The properties of $f(x)$ described in
Theorem~\ref{harr} will be discussed in the proof of
Lemma~\ref{sublimits-haar}. Since $f(x)$ is continuous and strictly
increasing in $(0,1)$,  the inverse function $f^{-1}$ of $f(x)$ is
well defined in the interval $(f(0), f(1))$.  It is also easy to
verify that the distribution function $F^*(y)$ is continuous in
$\mathbb{R}$ with density function
\begin{equation}\label{ff*}
f^*(y)=\frac{1}{f'(f^{-1}(y))}I(f(0)<y<f(1)).
\end{equation}

When $\varepsilon_j=1$ for $1\le j\le m_n$,  $Z_{m_n}$ is a product
of $m_n$ independent truncated Haar unitary matrices. Jiang and
Qi~\cite{JiangQi2019} obtained some limiting empirical distributions
under conditions imposed on $n_j$'s.  In particular, when $m=1$, the
work related to the eigenvalues $\bz_1, \cdots, \bz_n$ can be found
in \.{Z}yczkowski and Sommers~\cite{Zski}, Petz and
R\'effy~\cite{Petz2005}, Dong {\it et al.}~\cite{Dong}, Diaconis and
Evans~\cite{Evans}, and Jiang~\cite{Jiang09, Jiang2010}, Gui and
Qi~\cite{GQ2018}, and Miao and Qi~\cite{MQ2021}.

Under some special conditions on $n_j$'s,  conditions in
Theorem~\ref{harr} can be easily verified.  We will list some
results from Theorem~\ref{harr} in the following two remarks.

\noindent\textbf{Remark 4.} Let $m\ge 1$ be a fixed integer.  Assume
$\{\varepsilon_j,~1\le j\le m\}$ are constants taking values $1$ or
$-1$, $\{n_j,~1\le j\le m\}$ are $m$ integers such that $n_j>n$ for
all $1\le j\le m$, and there exist constants $\alpha_j\in [0,1]$
such that
\[
\lim_{n\to\infty}\frac{n}{n_j}=\alpha_j,~~~1\le j\le m.
\]

\noindent (i) If $\alpha_1=\cdots=\alpha_m=1$, then
$\lim_{n\to\infty}\Delta_n=0$. In this case,  by setting
$\gamma_n=2$ in Theorem~\ref{harr}(i), we have $\mu_n
\rightsquigarrow \mu=\mathrm{Unif[0,2\pi)}\otimes \delta_1$, which
implies
\[
\mu_n^*=\frac{1}{n}\sum^{n}_{j=1}\delta_{\bz_j/a_n^{1/2}}
\rightsquigarrow \mathrm{Unif}\{|z|=1\},
\]
where $a_n=\prod^m_{k=1}(\frac{n}{2n_k-n})^{\varepsilon_k}$.

\noindent (ii)  If $\min_{1\le j\le m}\alpha_j<1$, then
$\lim_{n\to\infty}\delta_{nj}=\frac1j\sum^m_{k=1}(-\varepsilon_k)^{j-1}\big(1-(\frac{\alpha_k}{2-\alpha_k})^j\big)$
for $j\ge 1$, and
$\lim_{n\to\infty}\Delta_n=\lim_{n\to\infty}\delta_{n1}=\sum^m_{k=1}\frac{2(1-\alpha_k)}{2-\alpha_k}\in
(0,\infty)$.   Now we can set $\gamma_n=2$ and condition
\eqref{limits} holds with
\[
\beta_j=\frac1{2j}\sum^m_{k=1}(-\varepsilon_k)^{j-1}\big(1-(\frac{\alpha_k}{2-\alpha_k})^j\big),
~~j\ge 1.
\]
It follows from Theorem~\ref{harr}(ii) that $\mu_n \rightsquigarrow
\mu=\mathrm{Unif[0,2\pi)}\otimes \nu$, where $\nu$ is a probability
measure with density function $f_v$ given by $f_v(y)=f^*(\log y)/y$,
$y>0$, where $f$ and $f^*$ are defined in \eqref{ff} and
\eqref{ff*}, respectively. Equivalently, we have
\[
\mu_n^*=\frac{1}{n}\sum^{n}_{j=1}\delta_{\bz_j/a_n^{1/2}}
\rightsquigarrow \mu^*
\]
with $a_n=\prod^m_{k=1}(\frac{n}{2n_k-n})^{\varepsilon_k}$, where
$\mu^*$ is a probability measure in $\mathbb{C}$ with density
function
\[
\frac{1}{2\pi|z|}f_v(|z|)=\frac{1}{2\pi|z|^2}f^*(\log |z|).
\]

\noindent\textbf{Remark 5.} Assume $m\to\infty$ as $n\to\infty$.
Then we have
\[
n\Delta_n=n\sum^m_{k=1}\frac{2(n_k-n)}{2(n_k-n)+n}\ge
n\sum^m_{k=1}\frac{2}{2+n}=\frac{2nm}{2+n},
\]
which implies $\lim_{n\to\infty}n\Delta_n=\infty$. Now we assume
\[
\lim_{n\to\infty}\max_{1\le j\le m}|\frac{n}{n_j}-\alpha|=0
\]
for some $\alpha\in [0,1)$, and
$\lim_{n\to\infty}\frac{p}{m}=\sigma\in [0,1]$, where $p=\#\{1\leq i
\leq m: \varepsilon_i=1\}$.  Now we set $\gamma_n=m$. It is easy to
verify \eqref{limits}. In fact, we have
\[
\beta_j=\lim_{n\to\infty}\frac{\delta_{nj}}{\gamma_n}=\lim_{n\to\infty}\frac{\delta_{nj}}{m}=
\left\{
  \begin{array}{ll}
   \frac{1}{j}\big(1-(\frac{\alpha}{2-\alpha})^{j}\big) , & \hbox{ if $j$ is odd}; \\
    \frac{1-2\sigma}{j}\big(1-(\frac{\alpha}{2-\alpha})^{j}\big), & \hbox{ if $j$ is even}.
  \end{array}
\right.
\]
The convergence of $\mu_n$ follows from Theorem~\ref{harr} (ii).


\medskip

{

\noindent\textbf{Remark 6.} The weak convergence in Theorems \ref{ginibre} and \ref{harr} imply that
\begin{eqnarray}\label{GFF}
\frac1n\sum^n_{j=1}\left[f(\mathbf{z}_j/a_n) - E_{\mu}f(\mathbf{z})\right] \to 0
\end{eqnarray}
almost surely for any bounded continuous function $f(z)$ defined on the complex plane, where $E_{\mu}f(\mathbf{z})$ is the mean value of $f(\mathbf{z})$ with random variable $\mathbf{z}$ following the limiting law $\mu.$ This is a kind of uniform law of large numbers. Notice that the eigenvalues of the random matrices form determinantal point processes with joint density functions in Propositions \ref{prop1} and \ref{prop2}, respectively. They have the same flavor of eigenvalues of matrices studied by Rider and Virag~\cite{Rider2007} in which Gaussian free fields are obtained as limits, precisely, the central limit theorems are derived for the object $\sum^n_{j=1}\left[f(\mathbf{z}_j/a_n) - E_{\mu}f(\mathbf{z})\right]$ for a class of smooth functions $f(\cdot).$ The two stories suggest Gaussian free fields could be established as limits if one tries to refine the law of large numbers aforementioned above  and studies the central limit theorems.

On the other hand, Webb~\cite{Webb2015} studies Gaussian multiplicative
chaos for Haar unitary matrices.  Berestycki~{\it et al.}~\cite{Berestycki2018} study the same phenomenon for random Hermitian matrices. Similar works are obtained by Hughes {\it et al.}~\cite{Hughes2001} and Fyodorov {\it et al.}~\cite{Fyodorov2016} Our product ensembles are different matrices from those but have similar features. It is conceivable that the analogues for our matrices should also hold. We leave it as a future work.
}

{
\noindent\textbf{Remark 7.} A recent breakthrough in quantum computing is showing the so-called quantum supreamcy by using photons; see, for instance, Yang {\it et al.}~\cite{Yang2020}, Yin {\it et al.}~\cite{Yin2020} and Zhong {\it et al.}~\cite{Zhong2020} The above research are based on the method ``Gaussian Boson sampling``, initially proposed
by Aaronson and Arkhipov~\cite{Aaronson2013}. Jiang and Ma~\cite{JiangMa2019} solved two conjectures by Aaronson and Arkhipov~\cite{Aaronson2013}; see also a comment
on p.~8 from Chabaud {\it et al.}~\cite{Chaud2021}  Jiang~\cite{Jiang2009b} also partially confirmed a conjecture
by Deshpande {\it et al.}~\cite{Deshpande2022} in a research on the Gaussian Boson sampling. Next let us briefly introduce the connection between this quantum algorithm and truncated Haar unitary matrices.

Unlike using the strings of $0$ and $1$ in classical quantum computing, quantum computing are realized through quantum gates, which are unitary matrices. Some common gates are, for example, Hadama gates, Pauli gates and Toffoli gates.  The tensors of different unitary matrices/gates have the same function as the strings of $0$ and $1$ in the classical quantum computing. Since there are many different gates/unitary matrices, an average target is an Haar unitary matrix. For a gate with big size, to understand the quantum advantage in the  the algorithm of the Gaussian Boson sampling, entries of truncation of Haar unitary matrices are needed to understood. This is the connection between the work by Jiang and Ma~\cite{JiangMa2019} and  Jiang~\cite{Jiang2009b} and the Gaussian Boson sampling. Although part of our work here is about the products of truncated Haar unitary matrices instead of tensor products, a  natural question is how and what we can understand about the tensor products of Haar unitary matrices or truncated Haar unitary matrices. For example, we may try to
obtain the limiting empirical  distributions of the  eigenvalues of the tensor products of truncations of Haar unitary matrices. Interesting enough, a recent work on quantum computing by Ma and Huang~\cite{Ma2025} shows that a specific quantum algorithm can create a Haar unitary matrix.  }

\section{Proofs}\label{proofs}

\subsection{Auxiliary lemmas}

In this subsection, we introduce some general results related to
eigenvalues that form a determinantal point process.

\begin{lemma}\label{nonlinear} (Theorem 1 in Jiang and Qi~\cite{JiangQi2019}).
 Let $\varphi(x)\geq 0$ be a measurable  function defined on $[0, \infty).$
Assume $(\bz_1, \cdots, \bz_n)\in \mathbb{C}^n$ has density $f(z_1,
\cdots, z_n)= C\cdot\prod_{1\leq j < k \leq n}|z_j-z_k|^2\cdot
\prod_{j=1}^n\varphi(|z_j|)$, where $C$ is a normalizing constant.
Let $Y_1, \cdots, Y_n$ be independent r.v.'s such that the density
of $Y_j$ is proportional to $y^{2j-1}\varphi(y)I(y\geq 0)$ for every
$1\leq j\leq n.$ Let $\mu_n$, $\nu_n$ and $\nu_n^*$ be defined as in
\eqref{mun}, \eqref{nun} and \eqref{mu*}, respectively. If $\{h_n\}$ are
measurable functions such that $\nu_n\rightsquigarrow\nu$ for some
probability measure $\nu$, then $\mu_n \rightsquigarrow \mu$ with
$\mu=\mathrm{Unif}[0, 2\pi]\otimes\nu$.
Taking $h_n(r)=r/a_n$, the conclusion still holds if ``$(\mu_n,
\nu_n, \mu, \nu)$" is replaced by  ``$(\mu_n^*, \nu_n^*, \mu^*,
\nu^*)$" where $\mu^*$ is the distribution of $Re^{i\Theta}$ with
$(\Theta, R)$ having the law of  $\mathrm{Unif}[0,
2\pi]\otimes\nu^*$.
\end{lemma}


\begin{lemma}\label{lemJQ} (Lemma 4 in Jiang and Qi~\cite{JiangQi2019})
Suppose $\{h_n(x);\, n\geq 1\}$ are measurable functions defined on
$[0,\infty)$ and $\nu_n$'s are defined as in \eqref{nun}. Let $Y_1,
\cdots, Y_n$ be as in Lemma~\ref{nonlinear} and $\nu$ be a
probability measure on $\mathbb{R}.$ Then $\nu_n\rightsquigarrow\nu$
if and only if
 \begin{equation}\label{v-limit}
\lim_{n\to\infty}\frac{1}{n}\sum^n_{j=1}\P(h_n(Y_j)\le r)=H(r)
\end{equation}
for every continuity point $r$ of $H(r)$, where $H(r):=\nu((-\infty,
r]),\, r \in \mathbb{R}$.
\end{lemma}

\begin{lemma}\label{ifandonlyif} Assume the conditions in
Lemma~\ref{nonlinear} hold and  $\{h_n\}$ are measurable functions
on $\mathbb{R}$. Then $\mu_n \rightsquigarrow \mu$ for some
probability measure $\mu$ if and only if \eqref{v-limit} holds for
some probability distribution $H$, and the relationship between
$\mu$ and $H$ is $\mu=\mathrm{Unif}[0, 2\pi]\otimes\nu$, where $\nu$
is a probability measure induced by $H$: $H(r)=\nu((-\infty, r]),\,
r \in \mathbb{R}$.
\end{lemma}

\noindent\textbf{Proof}. The sufficiency follows immediately from
Lemmas~\ref{nonlinear} and \ref{lemJQ}. Now we prove the necessity
part.  Recall that $\mu_n \rightsquigarrow \mu$ implies that for any
bounded and continuous function $u(\theta, x)$ defined on $\mathbb
R^2$
\[
\lim_{n\to\infty}\int
ud\mu_n=\lim_{n\to\infty}\frac1n\sum^n_{j=1}u(\Theta_j,
h_n(|\bz_j|))=\int ud\mu~~~\mbox{ with probability $1$}.
\]
In particular, for any bounded and continuous univariate function
$w(x)$, set $u(\theta, x)=w(x)$, the above equation implies
\[
\frac{1}{n}\sum^n_{j=1}w(h_n(|\bz_j|))\to \int wd\tau
\]
with probability $1$, where $\tau(A)=\mu(\mathbb{R}\times A)$ for
any measurable set $A$ in $\mathbb{R}$, which implies by using the
dominated convergence theorem
\[
\frac{1}{n}\sum^n_{j=1}\E w(h_n(|\bz_j|))\to \E\int wd\tau.
\]
From Lemma~1 in Jiang and Qi~\cite{JiangQi2019}, $\sum^n_{j=1}w(h_n(|\bz_j|))$
and $\sum^n_{j=1}w(h_n(Y_j))$ have the same distribution, we get
\eqref{v-limit}.  This completes the proof of
Lemma~\ref{ifandonlyif}. \eop

\subsection{Lemmas and proofs related to Theorem~\ref{ginibre}}

\begin{lemma}\label{t-moment} For the function $\varphi(z)$ given in
\eqref{dphi}, we consider  $\varphi(r),\, r> 0$, as its restriction on positive numbers.  Then,
\begin{equation}\label{moment2}
\int^\infty_0r^t\varphi(r)dr=\frac{\pi^{m-1}}{2}
\prod^m_{j=1}\Gamma\Big(\frac12\big(n+1+\varepsilon_j(t-n)\big)\Big),~~~~t>0.
\end{equation}
\end{lemma}

\noindent\textbf{Proof}. For $t>0$, we first compute $\int
|z|^{t-1}\varphi(|z|)dz$. By using the polar transformation we have
\begin{equation}\label{moment1}
 \int_{\mathbb{C}}
|z|^{t-1}\varphi(|z|)dz=2\pi\int^\infty_0r^t\varphi(r)dr.
\end{equation}
On the other hand, take $f(z)=|z|^{t-1}$ in \eqref{rep} and use \eqref{dphi} to see
\begin{eqnarray*}
\int_{\mathbb{C}}|z|^{t-1}\varphi(|z|)dz&=&\int_{\mathbb{C}^m}(|z_1|^{\varepsilon_1}\cdots|z_m|^{\varepsilon_m})^{t-1}
e^{-\sum^m_{j=1}|z_j|^2}\prod^m_{j=1}|z_j|^{(1-\varepsilon_j)(n-1)}\prod^m_{j=1}dz_j\\
&=&\int_{\mathbb{C}^m}
e^{-\sum^m_{j=1}|z_j|^2}\prod^m_{j=1}|z_j|^{n-1+\varepsilon_j(t-n)}\prod^m_{j=1}dz_j.
\end{eqnarray*}
Set $z_j=r_je^{i\theta_j}$ for $r_j\geq 0$ and $\theta_j \in [0, 2\pi)$. Then the above is equal to that
\begin{eqnarray*}
&&(2\pi)^m\int_{(0,\infty)^m}.e^{-\sum^m_{j=1}r_j^2}\prod^m_{j=1}r_j^{n-1+\varepsilon_j(t-n)}\prod^m_{j-1}r_j\prod^m_{j=1}dr_j\\
&=&(2\pi)^m\int_{(0,\infty)^m}
\prod^m_{j=1}\big(r_j^{n+\varepsilon_j(t-n)}e^{-r_j^2}\big)\prod^m_{j=1}dr_j.
\end{eqnarray*}
By changing the variables $s_j=r_j^2$ we have
\begin{eqnarray*}
\int_{\mathbb{C}}|z|^{t-1}\varphi(|z|)dz
 &=&\pi^m
\int_{(0,\infty)^m}\prod^m_{j=1}\big(s_j^{[n-1+\varepsilon_j(t-n)]/2}e^{-s_j}\big)\prod^m_{j=1}ds_j\\
&=&\pi^m\prod^m_{j=1}\Gamma\Big(\frac12\big(n+1+\varepsilon_j(t-n)\big)
\Big),
\end{eqnarray*}
which, together with \eqref{moment1}, yields \eqref{moment2}. \eop

\begin{lemma}\label{mgfoflogy} Let $Y_1, \dots, Y_n$ be the
independent random variables defined in Lemma~\ref{nonlinear} with $m, n, \epsilon_j$ and
$\varphi$ being given in \eqref{dphi}. Set
\begin{equation}\label{ajk}
\alpha_{j,k}=\frac12\big[n+1+\varepsilon_k(2j-1-n)\big], ~~~1\le
k\le m,~1\le j\le n.
\end{equation}
Then, for each $1\leq j \leq n$, the moment-generating function of $\log Y_j$ is
\begin{equation}\label{MGF}
M_{\log Y_j}(t)=
\frac{\prod^m_{k=1}\Gamma(\alpha_{j,k}+\frac12\varepsilon_kt)
\big)}{\prod^m_{k=1}\Gamma(\alpha_{j,k})},
~~~-2j<t<2(n+1-j).
\end{equation}
Assume $\{s_{j,k}, ~1\le k\le m,~ 1\le j\le n\}$ are independent
random variables such that $s_{j,k}$ has a Gamma\,$(\alpha_{j,k})$
distribution,   that is, $s_{j,k}$ has the density function
$y^{\alpha_{j,k}-1}e^{-y}I(y\ge 0)/\Gamma(\alpha_{j,k})$. Then $Y_j$
has the same distribution as
$\prod^m_{k=1}s_{j,k}^{\varepsilon_k/2}$ for each $1\le j\le n.$
\end{lemma}

\noindent\textbf{Proof}. Recall Lemma \ref{nonlinear}. The density
function of $Y_j$ is given by
\[
p_j(y)=\frac{y^{2j-1}\varphi(y)I(y\geq 0)}
{\int^\infty_0y^{2j-1}\varphi(y)dy}.
\]
Now the moment-generating function of $\log Y_j$ is $M_{\log
Y_j}(t)=\E(\exp(t\log Y_j))=\E(Y_j^t)$, which is the $t$-th moment of $Y_j$.
Therefore, we have from Lemma~\ref{t-moment} that
\bea\label{reirei}
M_{\log Y_j}(t)=\frac{\int^\infty_0y^{2j-1+t}\varphi(y)dy}
{\int^\infty_0y^{2j-1}\varphi(y)dy}=
\frac{\prod^m_{k=1}\Gamma(\frac12\big(n+1+\varepsilon_k(2j-1-n)+\varepsilon_kt)
)}{\prod^m_{k=1}\Gamma(\frac12(n+1+\varepsilon_k(2j-1-n))
)}
\eea
for  $-2j<t<2(n+1-j)$, which is obtained by considering $\varepsilon_k=\pm 1$ and
every variable in the Gamma function is positive. This  proves \eqref{MGF}.  Notice
\begin{eqnarray*}
  \E\exp\Big(\frac{1}{2}t\varepsilon_k\log s_{j,k}\Big)
  &=&\E\big(s_{j,k}^{t\varepsilon_k/2}\big)\\
&=&\frac{1}{\Gamma(\alpha_{j,k})}
  \int_0^{\infty}y^{\alpha_{j,k}+(1/2)t\varepsilon_k}e^{-y}\,dy\\
  &=& \frac{\Gamma\big(\alpha_{j,k}+(1/2)t\varepsilon_k\big)}{\Gamma(\alpha_{j,k})}.
\end{eqnarray*}
By setting $\xi=\frac12\sum^m_{k=1}\varepsilon_k\log s_{j,k}$, we know from independence and \eqref{ajk} that $\E e^{t\xi}$ is the same as the right hand side of \eqref{reirei} for  $-2j<t<2(n+1-j)$. By the uniqueness theorem on moment generating functions, we know $\log Y_j$ and $\xi$ have the same distribution.  Therefore,  $Y_j$ has the same distribution as
$\prod^m_{k=1}s_{j,k}^{\varepsilon_k/2}$.  \eop

\begin{lemma}\label{gamma-expansion}
For any $K>0$,
\begin{equation}\label{loggamma}
\log\frac{\Gamma(x+b)}{\Gamma(x)}=b\log x+O\Big(\frac{|b|}{x}\Big)
\end{equation}
uniformly over all $b\in [-K, K]$ as $x\to \infty$.
\end{lemma}

\noindent\textbf{Proof}. Let $\psi$ be the digamma function defined
by $\psi(x)=\frac{d}{dx}\log\Gamma(x)=\frac{\Gamma'(x)}{\Gamma(x)}$,
$x>0$.  It is known that
 \begin{equation}\label{digamma}
 \psi(t)=\log t-\frac{1}{2t}+O\Big(\frac{1}{t^2}\Big) ~~
\mbox{and} ~~\psi'(t)=\frac{1}{t}+O\Big(\frac1{t^2}\Big)
\end{equation}
as $t\to\infty$; see, e.g., Formulas 6.3.18 and 6.4.12 in Abramowitz and
Stegun~\cite{Abramowitz1972}. Then,
\[
\log\frac{\Gamma(x+b)}{\Gamma(x)}=\int^b_0\psi(x+t)dt=\int^b_0\log(x+t)dt+O\Big(\frac{|b|}{x}\Big).
\]
By expanding $\log(x+t)=\log x+\log(1+\frac{t}{x})=\log
x+O(\frac{t}{x})$ uniformly for $|t|\leq |b|$, we have
\[
\int^b_0\log(x+t)dt=b\log x+ O\Big(\frac{b^2}{x}\Big)=b\log x+
O\Big(\frac{|b|}{x}\Big).
\]
This proves \eqref{loggamma}.  \eop

\medskip


Let $Y_j$ be a random variables defined in Lemma~\ref{nonlinear}. We will study its mean and variance next. Easily,
\begin{equation}\label{ajk-case}
\alpha_{j,k}=(n+1)\frac{1-\varepsilon_k}{2}+j\varepsilon_k = \left\{
  \begin{array}{ll}
    j, & \hbox{ if }\varepsilon_k=1, \\
    n+1-j, & \hbox{ if } \varepsilon_k=-1.
  \end{array}
\right.
\end{equation}
Then we rewrite \eqref{MGF} as
\begin{equation}\label{helper}
M_{\log Y_j}(t)= \Big[\frac{\Gamma(j+t/2)}{\Gamma(j)}\Big]^p
\Big[\frac{\Gamma(n+1-j-t/2)}{\Gamma(n+1-j)}\Big]^{m-p}
\end{equation}
for  $-2j<t<2(n+1-j)$, where $p=\#\{i: \varepsilon_i=1, 1\le i\le
m\}$.

\begin{lemma}\label{v-of-Yj-ginibre}  Let $Y_1, \dots, Y_n$ be the
independent random variables defined in Lemma~\ref{nonlinear} with $m, n$ and
$\varphi$ being given in \eqref{dphi}. Assume $m=m_n$ depending on $n$. Set
$a_n=n^{2p-m}$. Then
\begin{equation}\label{mean}
\max_{n\delta-1\le j\le
n(1-\delta)+1}\left|a_n^{-1/m}\E(Y_j^{2/m})-G_{p/m}\big(jn^{-1}\big)\right|=O\Big(
\frac1n\Big)
\end{equation}
for any $\delta\in (0,\frac12)$, and
\begin{equation}\label{var}
\frac{1}{a_n^{2/m}}\max_{n\delta-1\le j\le
n(1-\delta)+1}\Var\big(Y_j^{2/m}\big)=O\Big(\frac1n\Big),
\end{equation}
which implies that
 \begin{equation}\label{2ndmoment}
\max_{n\delta-1\le j\le
n(1-\delta)+1}\E\left|a_n^{-1/m}Y_j^{2/m}-G_{p/m}\big(jn^{-1}\big)\right|=O\Big(
\frac1{\sqrt{n}}\Big).
\end{equation}
\end{lemma}

\noindent\textbf{Proof}. First, we have from \eqref{boundofG} that
\begin{equation}\label{kn-order}
\max_{n\delta-1\le j\le n(1-\delta)+1}G_{p/m}\Big(\frac{j}{n}\Big)=O(1)
\end{equation}
as $n\to\infty$. In what follows, we will assume $n\delta-1\le j\le n(1-\delta)+1$. Observe that $\E(Y_j^t)=M_{\log Y_j}(t)$. In view of
\eqref{helper} and Lemma~\ref{gamma-expansion} we have that, for any fixed
$t>0$,
\begin{equation}\label{ewuewo}
\log\E(Y_j^{t/m})=p\log\frac{\Gamma(j+(t/2m))}{\Gamma(j)}+(m-p)
\log \frac{\Gamma(n+1-j-(t/2m))}{\Gamma(n+1-j)}
\end{equation}
as $n$ is large enough such that $j+(t/2m)>0$ and $n+1-j-(t/2m)>0$. Recall the notation $a_n=n^{2p-m}$. Write $\log(n+1-j)=\log n+\log(1-\frac{j}{n})+\log (1+\frac{1}{n-j}).$
From Lemma \ref{gamma-expansion} and the formula $\log (1+\frac{1}{n-j})=O(\frac{1}{n-j})$ as $n\to\infty$, we know \eqref{ewuewo} is equal to
\begin{eqnarray*}
&&\frac{pt}{2m}\log
j-\frac{(m-p)t}{2m}\log(n+1-j)+O\Big(\frac{p}{mj}+\frac{m-p}{m(n-j)}\Big)\\
&=&\frac{tp}{2m}\log
\frac{j}{n}-\frac{(m-p)t}{2m}\log\Big(1-\frac{j}n\Big)+\frac{t}{2m}\log a_n+O\Big(\frac{1}{n}\Big)\\
&=&\frac{t}{2}\log
G_{p/m}\Big(\frac{j}{n}\Big)+\frac{t}{2m}\log a_n+O\Big(\frac1n\Big)
\end{eqnarray*}
uniformly over $n\delta-1\le j\le n(1-\delta)+1$ as $n\to\infty$,
which implies
\begin{equation}\label{gffdpo}
\max_{n\delta-1\le j\le
n(1-\delta)+1}\left|\log \frac{\E(Y_j^{t/m})}{a_n^{t/(2m)}}-\log G_{p/m}^{t/2}\big(jn^{-1}\big)\right|=
O\Big(\frac{1}{n}\Big)
\end{equation}
as $n\to\infty$. Set
\begin{eqnarray*}
u_{n,j}=a_n^{-t/(2m)}\E(Y_j^{t/m})\ \ \mbox{and}\ \
v_{n,j}=G_{p/m}^{t/2}\big(jn^{-1}\big).
\end{eqnarray*}
By \eqref{boundofG}, $0< v_{n,j}\leq C(t, \delta)$ uniformly for all
$n\delta-1\le j\le n(1-\delta)+1$, where $C(t, \delta)>0$ is a
finite constant depending on $t$ and $\delta.$ Evidently, the
identity $\log u_{n,j}=\log v_{n,j}+O(\frac{1}{n})$ implies that
$u_{n,j}=v_{n,j}\cdot e^{O(1/n)}=v_{n,j} + O(\frac{1}{n})$ uniformly
for all $n\delta-1\le j\le n(1-\delta)+1$. Consequently,
\eqref{gffdpo} concludes
\begin{eqnarray*}
\max_{n\delta-1\le j\le
n(1-\delta)+1}\left|a_n^{-t/(2m)}\E(Y_j^{t/m})-G_{p/m}^{t/2}\big(jn^{-1}\big)\right|=
O\Big(\frac{1}{n}\Big)
\end{eqnarray*}
as $n\to\infty$. Thus, \eqref{mean} follows by taking $t=2$. Now
choosing
 $t=4$, we have from \eqref{kn-order} that
 \begin{eqnarray*}
&&\max_{n\delta-1\le j\le
n(1-\delta)+1}a_n^{-2/m}\Var\big(Y_j^{2/m}\big)\\
&\le&\max_{n\delta-1\le j\le n(1-\delta)+1}
\left|a_n^{-2/m}\E(Y_j^{4/m})-G_{p/m}^2\big(jn^{-1}\big)\right|\\
&&~~~~~~~~~+\max_{n\delta-1\le j\le
n(1-\delta)+1}\left|\left[a_n^{-1/m}\E(Y_j^{2/m})\right]^2-G_{p/m}^2\big(jn^{-1}\big)\right|\\
&\le&O\big(n^{-1}\big)+\max_{n\delta-1\le j\le
n(1-\delta)+1}\left|a_n^{-1/m}\E(Y_j^{2/m})-G_{p/m}\big(jn^{-1}\big)\right|\cdot O(1)\\
&=&O\big(n^{-1}\big),
\end{eqnarray*}
proving \eqref{var}, where in the second inequality we use the fact $|a^2-b^2|\leq |a-b|\cdot (|a|+|b|)$ for any $a, b \in \mathbb{R}$, and \eqref{mean} and \eqref{kn-order}.  Finally,   by the Holder inequality, \eqref{mean} and
\eqref{var},
\begin{eqnarray*}
\E\left|a_n^{-1/m}Y_j^{2/m} -G_{p/m}\big(jn^{-1}\big)\right|
&\le&\Big[\E\Big(a_n^{-1/m}Y_j^{2/m} -G_{p/m}\big(jn^{-1}\big)\Big)^2\Big]^{1/2}\\
&=& \Big\{a_n^{-2/m}\Var\big(Y_j^{2/m}\big)+
\Big[a_n^{-1/m}\E\big(Y_j^{2/m}\big)-G_{p/m}\big(jn^{-1}\big)\Big]^2\Big\}^{1/2}.
\end{eqnarray*}
 We get \eqref{2ndmoment}. This completes the proof.\eop

Recall that a function $h$ is bounded and Lipschitz  on
$\mathbb{R}$ if,  for some $C>0$ and $K>0$,
\begin{equation}\label{Lip}
|h(x)|\le C~~\mbox{ and } |h(x)-h(y)|\le K|x-y|,  ~~~x, y\in
\mathbb{R}.
\end{equation}

\begin{lemma}\label{lemm-app} Let $Y_1, \dots, Y_n$ be the
independent random variables defined in Lemma~\ref{nonlinear} with $m, n$ and
$\varphi$ being given in \eqref{dphi}. Assume $m=m_n$ depending on $n$. Set
$a_n=n^{2p-m}$.  Then for any bounded and
Lipschitz function $h$ on $\mathbb{R}$ we have
\begin{equation}\label{app}
\lim_{n\to\infty}\Big|\frac{1}{n}\sum^n_{j=1}\E
h\left(\big(a_n^{-1}Y_j^2\big)^{1/m}\right)-\int^1_0h\big(G_{p/m}(x)\big)dx\Big|=0.
\end{equation}
\end{lemma}

\noindent\textbf{Proof}.  It follows from
Lemma~\ref{v-of-Yj-ginibre} that
\begin{eqnarray*}
\sup_{n\delta-1\le j\le n(1- \delta)+1}\E\Big|\big(a_n^{-1}Y_j^2\big)^{1/m}-G_{p/m}\big(\frac{j}{n}\big)\Big|=O\Big(\frac{1}{\sqrt{n}}\Big)
\end{eqnarray*}
as $n\to\infty$ for any given
$\delta\in (0,\frac{1}{2})$. Since $h$ is bounded and Lipschitz continuous on $\mathbb{R}$,
satisfying \eqref{Lip}, we have
 \begin{eqnarray}
&&\Big|\frac{1}{n}\sum^n_{j=1}\E
h\Big(\big(a_n^{-1}Y_j^2\big)^{1/m}\Big)-\int^1_0h\Big(G_{p/m}(x)\Big)dx\Big| \nonumber\\
&\le&\left|\frac{1}{n}\sum_{[n\delta]\le j\le n-[n\delta]} \E
h\left(\big(a_n^{-1}Y_j^2\big)^{1/m}\right)-\int^{1-\frac{[n\delta]-1}{n}}_{\frac{[n\delta]}{n}}
h\Big(G_{p/m}(x)\Big)dx\right|+4C\delta\nonumber\\
&\le&\left|\frac{1}{n}\sum_{[n\delta]\le j\le n-[n\delta]}\E\left[
h\Big(\big(a_n^{-1}Y_j^2\big)^{1/m}\Big)-h\Big(G_{p/m}\Big(\frac{j}{n}\Big)\Big)\right]\right|\nonumber\\
&&+\left|\int^{1-\frac{[n\delta]-1}{n}}_{\frac{[n\delta]}{n}}
\left[h\Big(G_{p/m}\Big(\frac{[nx]}{n}\Big)\Big)-h\Big(G_{p/m}(x)\Big)\right]dx\right|
+4C\delta,\label{dsjidso}
\end{eqnarray}
where the last step follows since
\begin{eqnarray*}
\frac{1}{n}\sum_{[n\delta]\le j\le n-[n\delta]}g\Big(\frac{j}{n}\Big)=\sum_{[n\delta]\le j\le n-[n\delta]}\int_{j/n}^{(j+1)/n}g\Big(\frac{[nx]}{n}\Big)\,dx
=\int_{[n\delta]/n}^{(n-[n\delta]+1)/n}g\Big(\frac{[nx]}{n}\Big)\,dx
\end{eqnarray*}
for any continuous function defined on $[0, 1]$. Easily, \eqref{dsjidso} bounded by
\begin{eqnarray*}
&&\frac{K}{n}\sum_{[n\delta]\le j\le
n-[n\delta]}\E\Big|(a_n^{-1}Y_j^2\big)^{1/m}-G_{p/m}\Big(\frac{j}{n}\Big)\Big|\\
&&
+K\int^{1-\frac{[n\delta]-1}{n}}_{\frac{[n\delta]}{n}}
\left|G_{p/m}\Big(\frac{[nx]}{n}\Big)-G_{p/m}(x)\right|dx+4K\delta\\
&\le&\frac{K}{n}\sum_{[n\delta]\le j\le
n-[n\delta]}O\Big(\frac{1}{\sqrt{n}}\Big)+K\cdot \sup_{\delta/2\le t\le
1-\delta/2}G_{p/m}'(t)\cdot \sup_{\delta/2\le x\le 1-
\delta/2}\Big|\frac{[nx]}{n}-x\Big|+4K\delta\\
&= &O\Big(\frac1{\sqrt{n}}\Big)+4K\delta
\end{eqnarray*}
by \eqref{boundofG} and the fact $\sup_{\delta/2\le x\le 1-
\delta/2}|\frac{[nx]}{n}-x|\leq \frac{1}{n}$. This and \eqref{dsjidso} conclude
\[
\limsup_{n\to\infty}\left|\frac{1}{n}\sum^n_{j=1}\E
h\Big(\big(a_n^{-1}Y_j^2\big)^{1/m}\Big)-\int^1_0h\Big(G_{p/m}(x)\Big)dx\right|
\le 4K\delta.
\]
Hence \eqref{app} follows by letting $\delta \downarrow 0$. The
proof is completed. \eop

\begin{lemma}\label{iff} Let $G^*_\alpha(x)$ be defined as in \eqref{G*}.
Assume $\{\alpha_n\}$ is  any sequence of numbers with $\alpha_n\in
[0,1]$. Then $G^*_{\alpha_n}$ converges weakly to a distribution $H$ if
and only if $\alpha_n$ has a limit, say $\alpha$, and $H$ must be
equal to $G^*_{\alpha}$.
\end{lemma}

\noindent\textbf{Proof}.  Let $\{F, ~F_n, n\ge 1\}$ be cumulative
distribution functions on $\mathbb{R}$.  It is known that $F_n$
converges weakly to $F$ if and only if
$\lim_{n\to\infty}\int^\infty_{-\infty}h(y)dF_n(y)=\int^\infty_{-\infty}h(y)dF(y)$
for all bounded Lipschitz continuous function $h$ on $\mathbb{R}$.
Observe that
\[
\int^\infty_{-\infty}h(y)dG^*_{\alpha}(y)=\int^1_0h(G_{\alpha}(x))dx
\]
for any bounded Lipschitz function $h$ and $\alpha\in [0,1]$.

Assume $\lim_{n\to\infty}\alpha_n=\alpha$. Then
$\lim_{n\to\infty}G_{\alpha_n}(x)=G_{\alpha}(x)$, $x\in (0,1)$.  By
using the dominated convergence theorem, we have
\[
\lim_{n\to\infty}\int^\infty_{-\infty}h(y)dG^*_{\alpha_n}(y)=\lim_{n\to\infty}\int^1_0h(G_{\alpha_n}(x))dx=\int^1_0h(G_{\alpha}(x))dx
=\int^\infty_{-\infty}h(y)dG^*_{\alpha}(y)
\]
for any bounded Lipschitz function $h$.  This proves the
sufficiency.

Now we assume that $G^*_{\alpha_n}$ converges weakly to a
distribution $H$, which implies
\[
\lim_{n\to\infty}\int^\infty_{-\infty}h(y)dG^*_{\alpha_n}(y)
=\int^\infty_{-\infty}h(y)dH(y).
\]
If the limit of $\alpha_n$ doesn't exist, we have two subsequential
limits, say $\lim_{n'\to\infty}\alpha_{n'}=\alpha\ne
\beta=\lim_{n''\to\infty}\alpha_{n''}$.  From the sufficiency, we
have
\[
\int^\infty_{-\infty}h(y)dH(y)=\int^\infty_{-\infty}h(y)dG^*_{\alpha}(y)=\int^\infty_{-\infty}h(y)dG^*_{\beta}(y)
\]
for all bounded Lipschitz functions $h$. By the Portmanteau theorem, bounded Lipschitz functions identify distributions. This implies
$G^*_{\alpha}(y)$ and $G^*_{\beta}(y)$ are the same distribution
functions, contradictory to the assumption $\alpha\ne \beta$.  This
completes the proof. \eop

\begin{lemma}\label{lem-sup} Let $Y_1, \dots, Y_n$ be the
independent random variables defined in Lemma~\ref{nonlinear} with $m, n$ and
$\varphi$ being given in \eqref{dphi}. Assume $m=m_n$ depending on $n$. Set
$a_n=n^{2p-m}$, $\alpha_n=p/m$, and
\begin{equation}\label{Hny}
H_n(y)=\frac{1}{n}\sum^n_{j=1}\P\Big(\big(a_n^{-1}Y_j^2\big)^{1/m}\le
y\Big), ~~~y\in \mathbb{R}.
\end{equation}
Then $\sup_{y \in \mathbb{R}}|H_n(y)-G^*_{\alpha_n}(y)| \to 0$ as
$n\to \infty$. Furthermore,  $\sup_{y \in
\mathbb{R}}|H_n(y)-G^*_{\alpha}(y)| \to 0$ if and only if
$\lim_{n\to\infty}\alpha_n=\alpha$.

\end{lemma}

\noindent\textbf{Proof}. Observe that
$\int^\infty_{-\infty}h(y)dG^*_{\alpha_{n}}(y)=\int_0^1h(G_{\alpha_n}(x))\,dx$
and
\begin{eqnarray*}
\int^\infty_{-\infty}h(y)dH_n(y)=\frac{1}{n}\sum^n_{j=1}\E
h\left(\big(a_n^{-1}Y_j^2\big)^{1/m}\right).
\end{eqnarray*}
By Lemma \ref{lemm-app},
\begin{eqnarray*}
\left|\int^\infty_{-\infty}h(y)dH_n(y)-\int^\infty_{-\infty}h(y)dG^*_{\alpha_{n}}(y)\right|\to 0
\end{eqnarray*}
as $n\to\infty.$ If $\lim_{n\to\infty}\alpha_n=\alpha$, then we have
from Lemma \ref{iff} that $G^*_{\alpha_n}$ converges weakly to
$G^*_{\alpha}$. As a consequence,
$\int^\infty_{-\infty}h(y)dG^*_{\alpha_{n}}(y)\to
\int^\infty_{-\infty}h(y)dG^*_{\alpha}(y).$ The two facts imply that
$\int^\infty_{-\infty}h(y)dH_n(y)\to
\int^\infty_{-\infty}h(y)dG^*_{\alpha}(y)$ for any bounded Lipschitz
function $h$. Thus, $H_n\to G^*_{\alpha}$ weakly. Since $H_n$ and
$G^*_{\alpha}$ are probability distribution functions and
$G^*_{\alpha}$ is continuous on $\mathbb{R}$, then $\sup_{y \in
\mathbb{R}}|H_n(y)-G^*_{\alpha}(y)| \to 0$; see, e.g., page 260 in Chow and Teicher~\cite{ChowandTeicher1978}. We
prove the sufficiency for  the  second conclusion of the lemma.  The
necessity simply follows from the sufficiency, otherwise if there
exists a subsequence of $n$, say $\{n'\}$ along which
$\alpha_{n'}\to\beta$, where $\beta\ne\alpha$, we have $H_{n'}$
converges to $G^*_{\beta}$ uniformly. This will lead to a
contradiction since $G^*_{\alpha}(y)\ne G^*_{\beta}(y)$ for some
$y$.

Now let us show the first claim in the lemma. In fact, suppose
$\sup_{y \in \mathbb{R}}|H_n(y)-G^*_{\alpha_n}(y)|$ does not go to
zero as $n\to \infty$. Then there is $\epsilon_0>0$ and  a
subsequence $\{n'\}$ such that
\begin{equation}\label{fewiufwei}
\sup_{y \in \mathbb{R}}|H_{n'}(y)-G^*_{\alpha_{n'}}(y)| \geq \epsilon_0.
\end{equation}
By definition, $\alpha_{n'} \in [0, 1]$. Then there is a further
subsequence $n''$ such that $\alpha_{n''}\to \alpha$. By the proved
conclusion, $\sup_{y \in \mathbb{R}}|H_{n''}(y)-G^*_{\alpha}(y)| \to
0$. Also, by Lemma \ref{iff}, $G^*_{\alpha_{n''}}$ converges weakly
to  $G^*_{\alpha}$. Since they all are continuous probability
distribution functions, we have $\sup_{y \in
\mathbb{R}}|G^*_{\alpha_{n''}}(y)-G^*_{\alpha}(y)| \to 0$ by the
Dini theorem. We have from the triangle inequality that $\sup_{y \in
\mathbb{R}}|H_{n''}(y)-G^*_{\alpha_{n''}}(y)|\to 0$. This
contradicts \eqref{fewiufwei}. Therefore, the first conclusion of
the lemmas holds. \eop


\noindent\textbf{Proof of Theorem~\ref{ginibre}}. Review the definition of $H_n(y)$ in \eqref{Hny}.  Since
\[
\frac{1}{n}\sum^n_{j=1}\P\Big(\big(a_n^{-1}Y_j^2\big)^{1/\gamma_n}\le
y\Big)=H_n(y^{\gamma_n/m_n}), ~~~y>0
\]
we have from the first claim in Lemma~\ref{lem-sup} that
\begin{equation}\label{appr}
\sup_y\Big|\frac{1}{n}\sum^{n}_{j=1}\P\Big(\big(a_n^{-1}Y_j^2\big)^{1/\gamma_{n}}\le
y\Big)-G^*_{\alpha_n}(y^{\gamma_{n}/m_{n}})\Big|\to 0
\end{equation}
 as $n\to\infty$, where $\alpha_n=p/m$. Now we prove parts (i) and (ii) below.

Note that
$h_n(r)=(a_n^{-1}r^2)^{1/\gamma_n}$.
 From Lemmas~\ref{nonlinear} and \ref{lemJQ}, we need to
show that
\begin{equation}\label{H}
\lim_{n\to\infty}\frac{1}{n}\sum^n_{j=1}\P\Big(\big(a_n^{-1}Y_j^2\big)^{1/\gamma_n}\le
y\Big)=H(y):=\left\{
     \begin{array}{ll}
       0, & \hbox{ if } y\in (0,1) \\
       1, & \hbox{ if } y>1
     \end{array}
   \right.
\end{equation}
since the probability distribution $H$ defined on the right-hand
above induces the probability measure $\delta_1$.   In view of
\eqref{appr},  it suffices to show that
\begin{equation}\label{sublimit}
\lim_{n\to\infty}G^*_{\alpha_n}(y^{\gamma_{n}/m_{n}})=H(y)=\left\{
     \begin{array}{ll}
       0, & \hbox{ if } y\in (0,1), \\
       1, & \hbox{ if } y>1
     \end{array}
   \right.
\end{equation}
under assumption that $m_n/\gamma_n\to 0$.  For any $y\in (0,1)$,
$y^{\gamma_n/m_n}\downarrow 0$ as $n\to\infty$, then it follows from
\eqref{OrderofG*} and \eqref{propertyof G*} that
\[
0\le G^*_{\alpha_n}(y^{\gamma_{n}/m_{n}})\le
G^*_1(y^{\gamma_{n}/m_{n}})\to 0.
\]
Similarly, for any $y>1$, $y^{\gamma_n/m_n}\uparrow \infty$ as
$n\to\infty$. We have from \eqref{OrderofG*} and \eqref{propertyof
G*} that
\[
1\ge G^*_{\alpha_n}(y^{\gamma_{n}/m_{n}})\ge
G^*_0(y^{\gamma_{n}/m_{n}})\to 1,
\]
proving \eqref{sublimit}.  This completes the proof of part (i).


Now we start to prove part (ii). Assume \eqref{twolimits} holds. We
follow the same lines as in the proof of Part (i), that is,
we need to show \eqref{H} with $H(y)=G^*_\alpha(y^{1/\beta})$,
$y>0$, or equivalently, we show
\[
\lim_{n\to\infty}G^*_{\alpha}(y^{\gamma_{n}/m_{n}})=G^*_\alpha(y^{1/\beta}),
~~y>0
\]
by using the second claim in Lemma~\ref{lem-sup} since
$\alpha_n\to\alpha$.   Because $\gamma_n/m_n\to 1/\beta$,
$y^{\gamma_n/m_n}\to y^{1/\beta}$ as $n\to\infty$. The above
equation is true since $G^*_{\alpha}$ is a continuous distribution
function.

Now we prove the necessity.  Assume $\mu_n \rightsquigarrow \mu$. In
view of Lemma~\ref{ifandonlyif},  there exist a probability
distribution function $H$ such that
\begin{equation}\label{apprH}
\lim_{n\to\infty}\frac{1}{n}\sum^n_{j=1}\P\Big(\big(a_n^{-1}Y_j^2\big)^{1/\gamma_n}\le
y\Big)=H(y)
\end{equation}
for all continuity points $y$ of $H$.  Here $H$ is a non-degenerate
distribution, otherwise the measure $\mu$ is singular since
$\mu=\mathrm{Unif}[0, 2\pi]\otimes \delta_a$ for some $a\ge 0$ from
Lemma~\ref{lemJQ}.  By combining \eqref{appr} and the above
equation, we have as $n\to\infty$
\begin{equation}\label{G*H}
G^*_{\alpha_n}(y^{\gamma_{n}/m_{n}})\to H(y)
\end{equation}
for all continuity points $y$ of $H$.

First, we see that $\beta_n:=m_n/\gamma_n$ is bounded away from
zero, otherwise there exists a subsequence along which $\beta_n$
converges to zero. From Part (i), we can conclude that the
subsequential limit of the left-hand side in \eqref{apprH} is a
degenerate distribution,  contradictory to the fact that $H$ is
non-degenerate.

Now we claim that $\beta_n$ is bounded from above. Otherwise, there
exists a subsequence of $n$ along which $\beta_n\to\infty$. We can
select its further subsequence, say $\{n'\}$, such that
$\alpha_{n'}\to \alpha\in [0,1]$.  Notice that
$\beta_{n'}\to\infty$. From Lemma~\ref{iff}, we have
$G^*_{\alpha_{n'}}$ converges weakly to $G^*_{\alpha}$ which is a
continuous distribution. Therefore,
\[
\sup_{y\ge
0}\big|G^*_{\alpha_{n'}}(y^{1/\beta_{n'}})-G^*_{\alpha}(y^{1/\beta_{n'}})\big|=\sup_{y\ge
0}\big|G^*_{\alpha_{n'}}(y)-G^*_{\alpha}(y)\big|\to 0
\]
as $n'\to\infty$, which coupled with \eqref{G*H} implies
\[
G^*_{\alpha}(y^{1/\beta_{n'}})\to H(y)
\]
for all continuity points $y$ of $H$. Since $\beta_{n'}\to\infty$,
$y^{1/\beta_{n'}}\to 1$ as $n'\to\infty$ for any $y>0$.  It follows
from the continuity of $G^*_{\alpha}$ that
$\lim_{n'\to\infty}G^*_{\alpha}(y^{1/\beta_{n'}})=G^*_{\alpha}(1)$
for $y>0$. This implies $H(y)=G^*_{\alpha}(1)$ for $y>0$.  Since $H$
is a probability distribution function, $\lim_{y\to\infty}H(y)=1$,
we have $H(y)=1$, $y>0$.  Therefore, we conclude that $H$ is a
probability distribution function degenerate at $0$, which yields a
contradiction.

Now we conclude that there exist two positive constants, $c_1$ and
$c_2$ such that $c_1\le \beta_n=m_n/\gamma_n\le c_2$ for all $n$ and
show \eqref{twolimits}.  We outline the proof. Consider $(\alpha_n,
\beta_n)$ as a vector. Then $(\alpha_n, \beta_n)\in [0,1]\times
[c_1,c_2]$. If \eqref{twolimits} is not true, then there exist two
subsequences of $\{n\}$ such that the two subsequential limits of
$(\alpha_n, \beta_n)$ are equal to $(a_1, b_1)$ and $(a_2, b_2)$,
respectively with $(a_1, b_1)\ne (a_2, b_2)$. Note that $b_1>0$ and
$b_2>0$. Along the two subsequences of $\{n\}$, we can use
\eqref{apprH} and \eqref{G*H} and  obtain that
$H(y)=G^*_{a_1}(y^{1/b_1})$ and $H(y)=G^*_{a_2}(y^{1/b_2})$ for all
$y>0$, which implies $G_{a_1}^{b_1}(x)=G_{a_2}^{b_2}(x)$ for $x\in
(0,1)$, that is
\[
x^{a_1b_1}(1-x)^{(a_1-1)b_1}=x^{a_2b_2}(1-x)^{(a_2-1)b_2}, ~~0<x<1.
\]
Therefore, we have
\[
a_1b_1=a_2b_2, ~~~(a_1-1)b_1=(a_2-1)b_2.
\]
Adding up two equations above leads to $b_1=b_2$ and $a_1=a_2$,
contradictory to the assumption that $(a_1, b_1)\ne (a_2, b_2)$.
The proof of part (ii) is completed.  \eop

\subsection{Lemmas and proofs related to Theorem~\ref{harr}}

Review the beta function $\mathrm{B}(a,b)$ defined in \eqref{beta}.
For any $a>0$ and $b>0$,  let $\mbox{Beta}(a,b)$ denote the beta distribution with density function
 \[
\frac{1}{\mathrm{B}(a,b)}t^{a-1}(1-t)^{b-1}I(0<t<1).
\]

\begin{lemma}\label{t-moment-Haar} Let the function $\varphi$ be given in
\eqref{dphi2}. Then
\begin{equation}\label{moment3}
\int^\infty_0r^t\varphi(r)dr=\frac{\pi^{m-1}}{2}
\prod^m_{j=1}\mathrm{B}\Big(\frac12\big(n+1+\varepsilon_j(t-n)\big),
n_j-n\Big)
\end{equation}
for any $t>0$.
\end{lemma}

\noindent\textbf{Proof}.  For $t>0$, we compute $\int
|z|^{t-1}\varphi(|z|)dz$ and then apply \eqref{moment1} to get
\eqref{moment3}.  By use of \eqref{dphi2} and \eqref{rep} and the
polar transformation
\begin{eqnarray*}
&&\int_{\mathbb{C}}|z|^{t-1}\varphi(|z|)dz\\
&=&\int_{\mathbb{C}^m}(|z_1|^{\varepsilon_1}\cdots|z_m|^{\varepsilon_m})^{t-1}
\prod^m_{j=1}(1-|z_j|^2)^{n_j-n-1}|z_j|^{(n-1)
(1-\varepsilon_j)}I(|z_j|\le 1)\prod^m_{j=1}dz_j\\
 &=&\int_{\mathbb{C}^m}
\prod^m_{j=1}(1-|z_j|^2)^{n_j-n-1} |z_j|^{n-1+\varepsilon_j(t-n)}I(|z_j|\le 1)\prod^m_{j=1}dz_j\\
&=&(2\pi)^m\int_{(0,1)^m}\prod^m_{j=1}(1-r_j^2)^{n_j-n-1}r_j^{n-1+\varepsilon_j(t-n)}\prod^m_{j-1}r_j\prod^m_{j=1}dr_j\\
&=&(2\pi)^m\int_{(0,1)^m}\prod^m_{j=1}r_j^{n+\varepsilon_j(t-n)}(1-r_j^2)^{n_j-n-1}\prod^m_{j=1}dr_j.
\end{eqnarray*}
By changing the variables $s_j=r_j^2$ we have
\begin{eqnarray*}
\int_{\mathbb{C}}|z|^{t-1}\varphi(|z|)dz
 &=&\pi^m
\int_{(0,\infty)^m}\prod^m_{j=1}\big(s_j^{(n-1+\varepsilon_j(t-n))/2}(1-s_j)^{n_j-n-1}\big)\prod^m_{j=1}ds_j\\
&=&\pi^m\prod^m_{j=1}\mathrm{B}\Big(\frac12\big(n+1+\varepsilon_j(t-n)\big),
n_j-n\Big),
\end{eqnarray*}
which, together with \eqref{moment1}, yields \eqref{moment3}. \eop

\begin{lemma}\label{rep-haar}
Let $Y_1, \dots, Y_n$ be the independent random variables defined in
Lemma~\ref{nonlinear} with $\varphi$ being given in \eqref{dphi2}.
Then the moment-generating function of $\log Y_j$ is
\begin{equation}\label{MGF2}
M_{\log Y_j}(t)=
\prod^m_{k=1}\frac{\mathrm{B}\big(\alpha_{j,k}+\frac{\varepsilon_k}{2}t,
n_k-n)}{\mathrm{B}\big(\alpha_{j,k}, n_k-n))} , ~~~-2j<t<2(n+1-j),
\end{equation}
where $\alpha_{j,k}$ is defined as in \eqref{ajk}. Assume
$\{s_{j,k}, ~1\le k\le m,~ 1\le j\le n\}$ are independent random
variables such that $s_{j,k}$ has a Beta($\alpha_{j,k}$, $n_k-n$)
distribution, that is, $s_{j,k}$ has a density function
$y^{\alpha_{j,k}-1}(1-y)^{n_k-n-1}I(0<y<1)/\mathrm{B}(\alpha_{j,k},
n_k-n)$. Then $Y_j$ has the same distribution as
$\prod^m_{k=1}s_{j,k}^{\varepsilon_k/2}$.
\end{lemma}

\noindent\textbf{Proof}. The proof is similar to that of
Lemma~\ref{mgfoflogy}, and the detail is omitted here. \eop

\begin{lemma}\label{v-of-Yj-unitrary}
Review that $Y_j$'s are defined in Lemma~\ref{rep-haar}.
Assume $m=m_n$ is any given sequence of positive integers that depend on
$n$. Set $a_n=\prod^m_{k=1}(\frac{n}{2n_k-n})^{\varepsilon_k}$ and
$\Delta_n=\sum^m_{k=1}\frac{2(n_k-n)}{2(n_k-n)+n}$. Then
\[
\E(\log Y_j)-\frac12\log
a_n=\frac12g_n(\frac{j}{n})+O(\frac{\Delta_n}{n})
\]
and
\[
\Var(\log Y_j)=O(\frac{\Delta_n}{n})
\]
uniformly for $n\delta\le j\le n(1-\delta)$ as $n\to\infty$, where
$\delta\in (0,\frac12)$ is any given number, and
\begin{equation}\label{gn}
g_n(x)=\sum^m_{k=1}\varepsilon_k\log\Big(1+2\varepsilon_k(x-\frac{1}{2})\Big)-\sum^m_{k=1}\varepsilon_k
\log\Big(1+\frac{2n\varepsilon_k}{2n_k-n}(x-\frac12)\Big), ~~x\in
(0,1).
\end{equation}
\end{lemma}

\noindent\textbf{Proof}. By using the formula
\[
\mathrm{B}(a,b)=\frac{\Gamma(a)\Gamma(b)}{\Gamma(a+b)},
\]
we have from \eqref{MGF2} and \eqref{ajk-case} that
\begin{eqnarray}\label{MGF-simple}
M_{\log Y_j}(t)&=&
\prod^m_{k=1}\frac{\Gamma\big(\alpha_{j,k}+\frac{\varepsilon_k}{2}t\big)}{\Gamma\big(\alpha_{j,k}\big)}
\prod^m_{k=1}\frac{\Gamma\big(\alpha_{j,k}+n_k-n\big)}
{\Gamma\big(\alpha_{j,k}+\frac{\varepsilon_k}{2}t+n_k-n\big)}\nonumber\\
&=&\prod_{k\in D^+}\Big(
 \frac{\Gamma(j+t/2)}{\Gamma(j)}\frac{\Gamma(n_k-n+j)}{\Gamma(n_k-n+j+t/2)}\Big)\\
&&\times\prod_{k\in
D^-}\Big(\frac{\Gamma(n+1-j-t/2)}{\Gamma(n+1-j)}\frac{\Gamma(n_k+1-j)}{\Gamma(n_k+1-j-t/2)}\Big)\nonumber
\end{eqnarray}
for $-2j<t<2(n+1-j)$, where \[ D^+=\{k: ~\varepsilon_k=1,~ 1\le k\le
m\}, ~~~D^-=\{k: ~\varepsilon_k=-1,~ 1\le k\le m\}.
\]

The first two derivatives of $M_{\log Y_j}(t)$ can be obtained as
follows
\[
M'_{\log Y_j}(t)=M_{\log Y_j}(t)\frac{d}{dt}\log M_{\log Y_j}(t),
\]
and
\begin{eqnarray*}
M''_{\log Y_j}(t)&=&M'_{\log Y_j}(t)\frac{d}{dt}\log M_{\log
Y_j}(t)+M_{\log Y_j}(t)\frac{d^2}{dt^2}\log M_{\log Y_j}(t)\\
&=&M_{\log Y_j}(t)\big(\frac{d}{dt}\log M_{\log
Y_j}(t)\big)^2+M_{\log Y_j}(t)\frac{d^2}{dt^2}\log M_{\log Y_j}(t).
\end{eqnarray*}
Moreover, it follows from \eqref{MGF-simple} that
\[
\frac{d}{dt}\log M_{\log Y_j}(t)=\frac12\Big\{\sum_{k\in D^+}\big(
 \psi(j+\frac{t}2)-\psi(n_k-n+j+\frac{t}2)\big)-\sum_{k\in
D^-}\big(\psi(n+1-j-\frac{t}2)-\psi(n_k+1-j-\frac{t}2)\big)\Big\}
\]
and
\[
\frac{d^2}{dt^2}\log M_{\log Y_j}(t)=\frac14\Big\{\sum_{k\in
D^+}\big(
 \psi'(j+\frac{t}2)-\psi'(n_k-n+j+\frac{t}2)\big)+\sum_{k\in
D^-}\big(\psi'(n+1-j-\frac{t}2)-\psi'(n_k+1-j-\frac{t}2)\big)\Big\},
\]
where $\psi$ is the digamma function defined in the beginning of the proof for Lemma~\ref{gamma-expansion}. These equations together with \eqref{digamma} yield
 \begin{eqnarray*}
 \E(\log Y_j)&=&M'_{\log Y_j}(0)\\
&=&\frac12\Big\{\sum_{k\in D^+}\big(
 \psi(j)-\psi(n_k-n+j)\big)-\sum_{k\in
D^-}\big(\psi(n+1-j)-\psi(n_k+1-j)\big)\Big\}\\
&=&\frac12\sum_{k\in
D^+}\Big\{\log\frac{j}{n_k-n+j}-\frac{1}{2}(\frac{1}{j}-\frac{1}{n_k-n+j})+O(\frac{1}{n^2})\Big\}\\
&&-\frac12\sum_{k\in
D^-}\Big\{\log\frac{n+1-j}{n_k+1-j}-\frac12(\frac1{n+1-j}-\frac1{n_k+1-j})+O(\frac1{n^2})\Big\}\\
&=&\frac12\sum_{k\in D^+}\log\frac{j}{n_k-n+j}-\frac12\sum_{k\in
D^-}\log\frac{n+1-j}{n_k+1-j}+O(\frac{1}{n}\sum^m_{k=1}\frac{n_k-n}{n_k})+O(\frac{m}{n^2})\\
&=&\frac12\sum_{k\in D^+}\log\frac{j}{n_k-n+j}-\frac12\sum_{k\in
D^-}\log\frac{n+1-j}{n_k+1-j}+O(\frac{1}{n}\sum^m_{k=1}\frac{n_k-n}{n_k})
\end{eqnarray*}
uniformly over $n\delta\le j\le n(1-\delta)$ as $n\to\infty$. In the
last step, we use the estimate
\[
\sum^m_{k=1}\frac{n_k-n}{n_k}=\sum^m_{k=1}(1-\frac{n}{n_k})\ge
m(1-\frac{n}{n+1})=\frac{m}{n+1}\ge \frac{m}{2n}.
\]
Similarly, we have
\begin{eqnarray*}
\Var(\log Y_j)&=&M''_{\log Y_j}(0)-(M'_{\log Y_j}(0))^2\\
&=&\frac14\Big\{\sum_{k\in D^+}\big(
 \psi'(j)-\psi'(n_k-n+j)\big)+\sum_{k\in
D^-}\big(\psi'(n+1-j)-\psi'(n_k+1-j)\big)\Big\}\\
&=&O(\frac{1}{n}\sum^m_{k=1}\frac{n_k-n}{n_k})
\end{eqnarray*}
uniformly over $n\delta\le j\le n(1-\delta)$ as $n\to\infty$.

It follows from Taylor's expansion that
\begin{eqnarray*}
\log\frac{n+1-j}{n_k+1-j}-\log\frac{n-j}{n_k-j}&=&\log(1+\frac{1}{n-j})-\log(1+\frac{1}{n_k-j})\\
&=&\frac{1}{n-j}-\frac1{n_k-j}+O(\frac{1}{n^2})
\end{eqnarray*}
uniformly over $n\delta\le j\le n(1-\delta)$ and $1\le k\le m$ as
$n\to\infty$. We conclude that
\[
\sum_{k\in D^-}\log\frac{n+1-j}{n_k+1-j}-\sum_{k\in
D^-}\log\frac{n-j}{n_k-j}=O(\frac1n\sum^m_{k=1}\frac{n_k-n}{n_k})
\]
uniformly over $n\delta\le j\le n(1-\delta)$ as $n\to\infty$. By
using the following trivial inequalities
\[
\frac12\Delta_n=\sum^m_{k=1}\frac{n_k-n}{2(n_k-n)+n}\le
\sum^m_{k=1}\frac{n_k-n}{n_k}\le
\sum^m_{k=1}\frac{2(n_k-n)}{2(n_k-n)+n}=\Delta_n,
\]
we obtain
\[
\E(\log Y_j)=\frac12b_{n,j}+O(\frac{\Delta_n}{n})
\]
and
\[
\Var(\log Y_j)=O(\frac{\Delta_n}{n})
\]
uniformly over $n\delta\le j\le n(1-\delta)$ and $1\le k\le m$ as
$n\to\infty$, where
\[
b_{n,j}=\sum_{k\in D^+}\log\frac{j}{n_k-n+j}-\sum_{k\in
D^-}\log\frac{n-j}{n_k-j}, ~~1\le j<n.
\]
Define $b_{n,n}=b_{n,n-1}$. It is easy to verify that $b_{n,j}$ is
increasing in $1\le j\le n-1$. We can also verify that
\begin{eqnarray*}
b_{n,j}-\log a_n&=&\sum_{k\in
D^+}\Big\{\log\big(1+2(\frac{j}{n}-\frac{1}{2})\big)-\log\big(1+\frac{2n}{2n_k-n}(\frac{j}{n}-\frac12)\big)\Big\}\\
&&-\sum_{k\in
D^-}\Big\{\log\big(1-2(\frac{j}{n}-\frac{1}{2})\big)-\log\big(1-\frac{2n}{2n_k-n}(\frac{j}{n}-\frac12)\big)\Big\}\\
&=&\sum^m_{k=1}\varepsilon_k\log\big(1+2\varepsilon_k(\frac{j}{n}-\frac{1}{2})\big)-\sum^m_{k=1}\varepsilon_k
\log\big(1+\frac{2n\varepsilon_k}{2n_k-n}(\frac{j}{n}-\frac12)\big)\\
&=&g_n(\frac{j}{n}).
\end{eqnarray*}
This completes the proof of the lemma. \eop

{   We still need to introduce a few lemmas before we prove Theorem~\ref{harr}.  Our objective is to establish a relationship  between the convergence of the empirical distribution based on scaled $Y_j$'s and convergence of $\delta_{n,j}$ in Theorem~\ref{harr}, where $Y_1,\cdots, Y_n$ are random variables defined in  Lemma~\ref{rep-haar}.
In particular,  we will show in Lemmas~\ref{gn-gn'} and \ref{complex-gn} that along any subsequence of $\{n\}$ such that \eqref{limits} holds,
the corresponding subsequential limit for   properly normalized function $g_n(x)$ is equal to the analytic function $f(x)$ defined in \eqref{ff}, where $g_n$ is defined as  \eqref{gn}. Lemma~\ref{complexlimit} is used to show the monotonicity of the limiting function $f(x)$, and
Lemmas~\ref{approximation} and \ref{sublimits-haar} show that
along any subsequence of $\{n\}$ such that condition \eqref{limits} in  Theorem~\ref{harr} holds,
the empirical distribution based on  $Y_1,\cdots, Y_n$
converges in distribution to the cumulative distribution function $F^*$ as defined in \eqref{F*}, which is
the inverse of the analytic function $f$. }

\begin{lemma}\label{gn-gn'}  Review that $g_n(x)$ is defined in
\eqref{gn} and $\Delta_n=\sum^m_{k=1}\frac{2(n_k-n)}{2(n_k-n)+n}$.  Then $g_n(x)$ is analytic on $(0,1)$ and both
$g_n(x)/\Delta_n$ and $g_n'(x)/\Delta_n$ are uniformly bounded in
any compact subset of $(0,1)$. Precisely, we have for any $\delta\in
(0,\frac12)$
\begin{equation}\label{gnbounded}
\frac{|g_n(x)|}{\Delta_n}\le
\frac{1}{1-2\delta}~~\mbox{if}~~|x-\frac12|\le\delta
\end{equation}
 and
\begin{equation}\label{gn'bounded}
\frac{2}{(1+2\delta)^2}\le
\frac{g'_n(x)}{\Delta_n}\le\frac{2}{(1-2\delta)^2}~~\mbox{if}~~|x-\frac12|\le\delta.
\end{equation}
\end{lemma}

\noindent\textbf{Proof}. For any $x\in (0,1)$, we have
\[
|2\varepsilon_k(x-\frac{1}{2})|=|2x-1|<1 ~\mbox{ and }~
|\frac{2n\varepsilon_k}{2n_k-n}(x-\frac12)|<|2x-1|<1. \]
 By use of
Taylor's expansion
$\log(1+t)=\sum^\infty_{j=1}(-1)^{j-1}\frac{t^j}{j}$ for $|t|<1$, we
get from \eqref{gn} that
\begin{eqnarray}\label{gn-expand}
g_n(x)&=&\sum^m_{k=1}\sum^\infty_{j=1}(-1)^{j-1}\varepsilon_k^{j+1}\frac{2^j(x-\frac{1}{2})^j}{j}-
\sum^m_{k=1}\sum^\infty_{j=1}(-1)^{j-1}\varepsilon_k^{j+1}\big(\frac{n}{2n_k-n}\big)^j\frac{2^j(x-\frac{1}{2})^j}{j}\nonumber\\
&=&\sum^\infty_{j=1}(-1)^{j-1}\sum^m_{k=1}\varepsilon_k^{j+1}\Big(1-\big(\frac{n}{2n_k-n}\big)^j
\Big)\frac{2^j(x-\frac{1}{2})^j}{j}\nonumber\\
&=&\sum^\infty_{j=1}\delta_{n,j}2^j(x-\frac{1}{2})^j,
\end{eqnarray}
where $\delta_{n,j}$ is defined in \eqref{deltanj}.

Since $0\le 1-t^j=(1-t)\sum^{j-1}_{k=0}t^j\le j(1-t)$ for $t\in
[0,1]$,
\begin{equation}\label{boundedratio}
 |\delta_{n,j}|\le
\frac1j\sum^m_{k=1}\Big(1-\big(\frac{n}{2n_k-n}\big)^j\Big)\le\sum^n_{k=1}(1-\frac{n}{2n_k-n})=\Delta_n.
\end{equation}
Then it follows from \eqref{gn-expand} that for any $\delta\in
(0,\frac12)$
\[
\frac{|g_n(x)|}{\Delta_n}\le
\sum^\infty_{j=1}2^j|x-\frac{1}{2}|^j\le\sum^\infty_{j=1}(2\delta)^j\le
\frac{1}{1-2\delta}
\]
if  $|x-\frac12|\le \delta$, proving \eqref{gnbounded}.

Next, we have from \eqref{gn} that
\begin{eqnarray*}
g_n'(x)&=&\sum^m_{k=1}\frac{2}{1+2\varepsilon_k(x-\frac{1}{2})}-\sum^m_{k=1}\frac{\frac{2n}{2n_k-n}}{1+\frac{2n\varepsilon_k}{2n_k-n}(x-\frac12)}\\
&=&2\sum^m_{k=1}(1-\frac{n}{2n_k-n})\frac{1}{\big(1+2\varepsilon_k(x-\frac12)\big)\big(1+\frac{2n\varepsilon_k}{2n_k-n}(x-\frac12)\big)}.
\end{eqnarray*}
Then by using estimates
\[
(1-2\delta)^2\le
\big(1+2\varepsilon_k(x-\frac12)\big)\big(1+\frac{2n\varepsilon_k}{2n_k-n}(x-\frac12)\big)\le
(1+2\delta)^2 ~~\mbox{for} ~~|x-\frac{1}{2}|\le \delta,
\]
we obtain for any $|x-\frac12|\le \delta$
 \[
 \frac{2}{(1+2\delta)^2}\Delta_n=\frac{2}{(1+2\delta)^2}\sum^m_{j=1}(1-\frac{n}{2n_k-n})\le g_n'(x)\le
\frac{2}{(1-2\delta)^2}\Delta_n,
\]
proving \eqref{gn'bounded}.
 \eop

Now we extend the domain of functions $g_n(x)$ to the complex plane.
More generally, for any sequence of positive numbers $\{\gamma_n\}$,
we define
\[
f_n(z)=\frac{g_n(z)}{\gamma_n}=\sum^\infty_{j=1}c_{nj}\big(2(z-1)\big)^j,~~~|z-\frac12|<\frac12,
~z\in \mathbb{C},
\]
where $c_{nj}=\delta_{n,j}/\gamma_n$.  In view of
\eqref{boundedratio},
$|c_{nj}|=\frac{\Delta_n}{\gamma_n}\frac{|\delta_{n,j}|}{\Delta_n}\le
\frac{\Delta_n}{\gamma_n}$ for all $j\ge 1$. If
$\frac{\Delta_n}{\gamma_n}$ is bounded, then $f_n(z)$ is uniformly
bounded in the disk $\{z\in \mathbb{C}: |z-\frac{1}{2}|\le \delta\}$
for any $\delta\in (0,\frac{1}{2})$, and $f_n(z)$ is analytic in
$\{z\in\mathbb{C}: |z-\frac{1}{2}|<\frac12\}$.


\begin{lemma}\label{complex-gn} Let $\{\gamma_n\}$ be a sequence of positive
numbers. Assume $\{n_s\}$ is a given subsequence of $\{n\}$ such
that
\begin{equation}\label{sublimits}
\lim_{s\to\infty}c_{n_sj}=\lim_{s\to\infty}\frac{\delta_{n_s,j}}{\gamma_{n_s}}=:\beta_j,
~~~j\ge 1.
\end{equation}
Write $\boldsymbol{\beta}=(\beta_1, \beta_2, \cdots)$ and define
\[
f_{\boldsymbol{\beta}}(z)=\sum^{\infty}_{j=1}\beta_j\big(2(z-1)\big)^j,~~~|z-\frac12|<\frac12,~
z\in \mathbb{C}.
\]
 Then $f_{n_s}(z)$ converges to $f_{\boldsymbol{\beta}}(z)$ uniformly on $\{z\in
\mathbb{C}: |z-\frac{1}{2}|\le \delta\}$ for any $\delta\in
(0,\frac{1}{2})$.
\end{lemma}

\noindent\textbf{Proof}. Set $c=\sup_{s\ge 1}\frac{\Delta_{n_s}}{\gamma_{n_s}}$, which is finite from assumption \eqref{sublimits} with $j=1$. For any fixed integer $k\ge 2$, we have for
$|z-\frac12|\le \delta<\frac12$
\begin{eqnarray*}
|f_{n_s}(z)-f_{\boldsymbol{\beta}}(z)|&\le&
|\sum^k_{j=1}(c_{n_sj}-\beta_j)\big(2(z-1)\big)^j|+|\sum^\infty_{j=k+1}c_{n_sj}\big(2(z-1)\big)^j-
\sum^\infty_{j=k+1}\beta_j\big(2(z-1)\big)^j|\\
&\le&\max_{2\le j\le
k}|c_{n_sj}-\beta_j|+2c\sum^\infty_{j=k+1}(2\delta)^j\\
&\le&\max_{2\le j\le
k}|c_{n_sj}-\beta_j|+\frac{2c(2\delta)^{k+1}}{1-2\delta},
\end{eqnarray*}
which
leads to
\[
\limsup_{s\to\infty}\sup_{|z-\frac12|\le
\delta}|f_{n_s}(z)-f_{\boldsymbol{\beta}}(z)|\le
\frac{2c(2\delta)^{k+1}}{1-2\delta}.
\]
The right-hand side tends to zero if we let $k$ tend to infinity.
This completes the proof of the lemma. \eop

\begin{lemma}\label{complexlimit} (Theorem 10.28 in Rudin~\cite{RudinComplex}) Suppose $f_j$ is analytic on open set $\Omega\subset \mathbb{C}$ for $j=1,2,\cdots,$ and $f_j\rightarrow f$ uniformly on each compact subset of $\Omega$ as $j\to\infty$. Then $f$ is analytic on $\Omega$, and $f'_j\rightarrow f'$ uniformly on any compact subset on $\Omega$ as $j\to\infty$. \end{lemma}

\begin{lemma}\label{approximation}
Assume $m=m_n$ is any given sequence of positive integers that
depend on $n$. Set
$a_n=\prod^m_{k=1}(\frac{n}{2n_k-n})^{\varepsilon_k}$ and
$\Delta_n=\sum^m_{k=1}\frac{2(n_k-n)}{2(n_k-n)+n}$, and assume
$\gamma_n$ is a sequence of positive numbers such that
\begin{equation}\label{Delta-gamma}
\lim_{n\to\infty}\frac{\Delta_n}{n\gamma_n^2}=0~\mbox{and}~\lim_{n\to\infty}\frac{\Delta_n}{n\gamma_n}=0.
\end{equation}
 Then for any bounded Lipschitz continuous
function $h$ on $\mathbb{R}$ we have
\begin{equation}\label{limit}
\lim_{n\to\infty}\Big|\frac{1}{n}\sum^n_{j=1}\E
h\Big(\frac1{\gamma_n}\big(2\log Y_j-\log
a_n\big)\Big)-\int^1_0h\Big(g_n(x)/\gamma_n\Big)dx\Big|=0,
\end{equation}
where $Y_j$'s are independent random variables defined in  Lemma~\ref{rep-haar},  and function $g_n$ is defined in \eqref{gn} or equivalently in \eqref{gn-expand}.
\end{lemma}

\noindent\textbf{Proof}.  It follows from
Lemma~\ref{v-of-Yj-unitrary} that
\begin{eqnarray*}
\E\big|2\log Y_j-\log a_n-g_n(\frac{j}{n})\big|&\le&
2\Big\{\E\Big(\log Y_j-\frac12\log a_n-\frac12g_n(\frac{j}{n})\Big)^2\Big\}^{1/2}\\
&=&
\Big\{\Var(\log Y_j)+\Big(\E\log Y_j-\frac12\log a_n-\frac12g_n(\frac{j}{n})\Big)^2\Big\}^{1/2}\\
&=&O(\sqrt{\frac{\Delta_n}{n}}+\frac{\Delta_n}{n}).
\end{eqnarray*}
uniformly over $n\delta/2\le j\le n(1- \delta/2)$ for any given
$\delta\in (0,\frac{1}{2})$.  Therefore,  for bounded Lipschitz
continuous function $h$ satisfying \eqref{Lip}, we have from
Lemma~\ref{gn-gn'} that
 \begin{eqnarray*}
&&\Big|\frac{1}{n}\sum^n_{j=1}\E h\Big(\frac1{\gamma_n}\big(2\log
Y_j-\log
a_n\big)\Big)-\int^1_0h\Big(g_n(x)/\gamma_n\Big)dx\Big|\\
&\le&\Big|\frac{1}{n}\sum_{[n\delta]\le j\le n-[n\delta]} \E
h\Big(\frac1{\gamma_n}\big(2\log Y_j-\log
a_n\big)\Big)-\int^{1-\frac{[n\delta]-1}{n}}_{\frac{[n\delta]}{n}}h\Big(g_n(x)/\gamma_n\Big)dx\Big|+4C\delta\\
&\le&\Big|\frac{1}{n}\sum_{[n\delta]\le j\le n-[n\delta]}\E\Big\{
h\Big(\frac1{\gamma_n}\big(2\log Y_j-\log
a_n\big)\Big)-h\Big(g_n\big(\frac{j}{n}\big)/\gamma_n\Big)\Big\}\Big|\\
&&+\Big|\int^{1-\frac{[n\delta]-1}{n}}_{\frac{[n\delta]}{n}}\Big\{h\Big(g_n(\frac{[nx]}{n})/\gamma_n\Big)-h\Big(g_n(x)/\gamma_n\Big)\Big\}dx\Big|
+4C\delta\\
&\le&\frac{K}{\gamma_n}\frac{1}{n}\sum_{[n\delta]\le j\le
n-[n\delta]}\E\Big|2\log Y_j-\log
a_n-g_n(\frac{j}{n})\Big|\\
&&+\frac{K}{\gamma_n}\Big|\int^{1-\frac{[n\delta]-1}{n}}_{\frac{[n\delta]}{n}}\big|g_n(\frac{[nx]}{n})-g_n(x)\big|dx+4C\delta\\
&\le&\frac{K}{\gamma_n}\frac{1}{n}\sum_{[n\delta]\le j\le
n-[n\delta]}O(\sqrt{\frac{\Delta_n}{n}}+\frac{\Delta_n}{n})+\frac{K\Delta_n}{\gamma_n}\sup_{\delta/2\le
t\le 1-\delta/2}\frac{|g_n'(t)|}{\Delta_n}\sup_{\delta/2\le x\le 1-
\delta/2}|\frac{[nx]}{n}-x| +4C\delta\\
&=&O(\sqrt{\frac{\Delta_n}{n\gamma_n^2}}+\frac{\Delta_n}{n\gamma_n})+4C\delta,
\end{eqnarray*}
which implies
\[
\limsup_{n\to\infty}\Big|\frac{1}{n}\sum^n_{j=1}\E
h\Big(\frac1{\gamma_n}\big(2\log Y_j-\log
a_n\big)\Big)-\int^1_0h\Big(g_n(x)/\gamma_n\Big)dx\Big|  \le
4C\delta.
\]
which yields \eqref{limit} by letting $\delta$ tend to $0$. The
proof of the lemma is completed. \eop

\begin{lemma}\label{sublimits-haar}
Let $\{\gamma_n\}$ be a sequence of positive numbers and define
\begin{equation}\label{Fn}
F_n(y)=\frac{1}{n}\sum^n_{j=1}\P\Big(\frac{1}{\gamma_n}(2\log
Y_j-\log a_n\big)\le y\Big), ~~~y\in \mathbb{R}
\end{equation}
for random variables $Y_j$'s given
in Lemma~\ref{rep-haar}.
Assume $\lim_{n\to\infty}n\Delta_n=\infty$. Let $\{n_s\}$ be a
subsequence of $\{n\}$ such that \eqref{sublimits} holds with
$\beta_1\in (0,\infty)$. Then
\[
\lim_{s\to\infty}\sup_y|F_{n_s}(y)-F^*_{\boldsymbol\beta}(y)|=0,
\]
where $\boldsymbol{\beta}$ denotes the vector $(\beta_1, \beta_2,
\cdots)$ with $\beta_j=\lim_{s\to\infty}\delta_{n_sj}/\gamma_{n_s}$
for $j\ge 1$, $F^*_{\boldsymbol\beta}$ is defined in \eqref{F*} with
$f$ being replaced by $f_{\boldsymbol{\beta}}$ which is defined in Lemma~\ref{complex-gn}.
\end{lemma}

\noindent\textbf{Proof}.  We have from Lemma~\ref{complex-gn} that
\[
\lim_{s\to\infty}f_{n_s}(x)=f_{\boldsymbol{\beta}}(x), ~~ x\in
(0,1).
\]
Note that $\Delta_n=\delta_{n,1}$ and we assume
$n\Delta_n\to\infty$. Condition that $\delta_{n,1}/\gamma_{n}$ converges
to $\beta_1\in (0,\infty)$ along the subsequence $\{n_s\}$ implies
\[
\frac{\Delta_{n_s}}{n_s\gamma_{n_s}^2}=
\frac{1}{n_s\Delta_{n_s}}\Big(\frac{\delta_{n_s,1}}{\gamma_{n_s}}\Big)^2\to 0 ~\mbox{ and }~
\frac{\Delta_{n_s}}{n_s\gamma_{n_s}}=
\frac{1}{n_s}\frac{\delta_{n_s,1}}{\gamma_{n_s}}\to 0,
\]
that is,
\eqref{Delta-gamma} holds along the subsequence $\{n_s\}$.  Therefore, we can
apply the convergence \eqref{limit} in Lemma~\ref{approximation}
along the same subsequence. Therefore, for any bounded Lipschitz function
$h$ we have
\begin{equation}\label{integral}
\lim_{s\to\infty}\int^\infty_{-\infty}
h(y)dF_{n_s}(y)=\lim_{s\to\infty}\int^1_0 h(f_{n_s}(x))dx=\int^1_0
h(f_{\boldsymbol{\beta}}(x))dx
\end{equation}
in view of the dominated convergence theorem.   Using
Lemmas~\ref{complex-gn} and \ref{complexlimit}, we have
\[
f_{\boldsymbol{\beta}}'(x)=\lim_{s\to\infty}f_{n_s}'(x)=\lim_{s\to\infty}\frac{\Delta_{n_s}}{\gamma_{n_s}}\frac{f_{n_s}'(x)}{\Delta_{n_s}}
=\beta_1\lim_{s\to\infty}\frac{f_{n_s}'(x)}{\Delta_{n_s}},
\]
which together with Lemma~\ref{gn-gn'} implies
\[
\frac{2\beta_2}{(1+2\delta)^2}\le f_{\boldsymbol{\beta}}'(x)\le
\frac{2\beta_1}{(1-2\delta)^2}, ~~|x-\frac12|\le \delta
\]
for any $\delta\in (0,\frac12)$. Therefore,
$f_{\boldsymbol{\beta}}(x)$ is analytic and strictly increasing in
$(0,1)$.

Now define $f_{\boldsymbol{\beta}}(0)=\lim\limits_{x\downarrow 0}f_{\boldsymbol{\beta}}(x)$
and $f_{\boldsymbol{\beta}}(1)=\lim\limits_{x\uparrow 1}f_{\boldsymbol{\beta}}(x)$. Let $U$ be a random variable uniformly distributed over $(0,1)$,
Then the cumulative distribution function of $f_{\boldsymbol{\beta}}(U)$ is $F_{\boldsymbol{\beta}}^*(y)$ which is defined in \eqref{F*} with
$f$ being replaced by $f_{\boldsymbol{\beta}}$. Therefore, it follows from
\eqref{integral} that
\[
\lim_{s\to\infty}\int^\infty_{-\infty}
h(y)dF_{n_s}(y)=\E h\big(f_{\boldsymbol{\beta}}(U)\big)
=\int^{\infty}_{-\infty}
h(y)dF_{\boldsymbol{\beta}}^*(y),
\]
which implies  $F_{n_s}$ converges weakly to
$F_{\boldsymbol{\beta}}^*$. The lemma is obtained since
$F_{\boldsymbol{\beta}}^*$ is continuous in $\mathbb{R}$
as the cumulative distribution function of continuous random variable
$f_{\boldsymbol{\beta}}(U)$.
\eop

\noindent\textbf{Proof of Theorem~\ref{harr}}. Note that $F_n$  defined in \eqref{Fn} is a probability distribution
of a continuous random variable and it also depends on $\gamma_n$.
To reflect its dependence on  $\gamma_n$, we use $W_{n, \gamma_n}$
to denote a random variable with distribution function $F_n$.   Then
the distribution of $\overline{W}_{n,\lambda_n}:=\exp(W_{n,
\lambda_n})$ is given by
\[
F_{\overline{W}_{n,\gamma_n}}(y)=F_n(\log
y)=\frac{1}{n}\sum^{n}_{j=1}\P\Big((a_n^{-1}Y_j^2)^{1/\gamma_n}\le
y\Big), ~~~y>0.
\]

\noindent{\underline{Part (i)}.  From \eqref{gnbounded}, we
have
$\frac{g_n(x)}{\gamma_n}=\frac{\Delta_n}{\gamma_n}\frac{g_n(x)}{\Delta_n}\to
0$ for any $x\in (0,1)$. We note that
\[
\Delta_n=\sum^m_{k=1}\frac{2(n_k-n)}{2(n_k-n)+n}\ge
\sum^m_{k=1}\frac{2}{2+n}\ge \frac{m}{2n}.
\]
This together with assumption $\Delta_n/\gamma_n\to 0$ implies that
the conditions in Lemma~\ref{approximation} hold. Then for any
bounded Lipschitz function $h$, by using \eqref{limit} and the
dominated convergence theorem, we obtain
\[
\lim_{n\to\infty}\int^\infty_{0}h(y)dF_n(y)=h(0),
\]
which implies $W_{n,\lambda_n}$ converges in probability to zero, or
equivalently, $\overline{W}_{n,\lambda_n}=\exp(W_{n, \lambda_n})$
converges in probability to one. In view of Lemma~\ref{lemJQ}, this
proves  Part (\textbf{i}).

\noindent{\underline{Part (ii)}. We will first prove the
sufficiency, that is, $\mu_n \rightsquigarrow
\mu=\mathrm{Unif[0,2\pi)}\otimes \nu$ under condition
\eqref{limits}, where $\mu$ is non-degenerate probability measure
given in Theorem~\ref{harr}.


We have assumed that $\lim_{n\to\infty}n\Delta_n=\infty$. From
Lemma~\ref{sublimits-haar},
\[
\lim_{n\to\infty}F_{\overline{W}_{n,\gamma_n}}(y)=\lim_{n\to\infty}F_n(\log
y)=F^*_{\boldsymbol\beta}(\log y), ~~y>0.
\]
Note that $F^*_{\boldsymbol\beta}(\log y)$, $y>0$, is a continuous
probability distribution function for a positive random variable.
Let $\nu$ be a probability measure induced by
$F^*_{\boldsymbol\beta}(\log y)$. Then, in view of
Lemma~\ref{lemJQ}, we obtain that $\mu_n \rightsquigarrow
\mu=\mathrm{Unif[0,2\pi)}\otimes \nu$.

To prove the necessity, we will show
$\sup_n(\Delta_n/\gamma_n)<\infty$ first. If this is not true, then
there exists a subsequence, say $\{n_s\}$, of $\{n\}$ such that
$r_s=:\Delta_{n_s}/\gamma_{n_s}\to\infty$ as $s\to\infty$.  In view
of \eqref{boundedratio}, we see that $|\delta_{n_s,j}|/\Delta_{n_s}$
are bounded by $1$.  Using the diagonal argument, there exists a
subsequence of $\{n_s\}$ along which $\delta_{n_s,j}/\Delta_{n_s}$
has a limit for each $j\ge 1$. For brevity, we simply assume the
limit of $\delta_{n_s,j}/\Delta_{n_s}$ exists for each $j\ge 1$ as
$s\to\infty$, and denote these limits by $\beta_j$, $j=1,2, \cdots$
with $\beta_1=1$.  Set $\underline{\boldsymbol{\beta}}=(1, \beta_2,
\beta_3,\cdots)$. Now we apply Lemma~\ref{sublimits-haar} with
$\gamma_n=\Delta_n$ and have
\begin{equation}\label{beta1=1}
\lim_{s\to\infty}\sup_y|F_{W_{n_s,\Delta_{n_s}}}(y)-F^*_{\underline{\boldsymbol\beta}}(y)|=0.
\end{equation}
Note that $f_{\underline{\boldsymbol\beta}}(x)$ is analytic and
strictly increasing in  $(0,1)$, and
$f_{\underline{\boldsymbol\beta}}(\frac12)=0$, we have
$F^*_{\underline{\boldsymbol\beta}}(0)=\frac{1}{2}$.

Since we have assumed that $\nu$ is a non-singular probability
measure,  $\nu$ is a non-singular probability measure, which implies
by using Lemma~\ref{lemJQ} that $\overline{W}_{n_s,\gamma_{n_s}}$
converges in distribution to a non-degenerate distribution. Note
that a degenerate distribution induces a singular probability
measure.

We notice that,  for $y>0$,
\begin{eqnarray*}
F_{\overline{W}_{n_s,\gamma_{n_s}}}(y)
&=&\frac{1}{n_s}\sum^{n_s}_{j=1}\P\Big(\big(a_{n_s}^{-1}Y_j^2\big)^{1/\gamma_{n_s}}\le
y\Big)\\
&=&\frac{1}{n_s}\sum^{n_s}_{j=1}\P\Big(\frac{1}{\Delta_{n_s}}\big(2\log
Y_j-\log a_{n_s}\big)\le
\frac{\gamma_{n_s}}{\Delta_{n_s}}\log y\Big)\\
&=&F_{W_{n_s,\Delta_{n_s}}}\big(\frac{1}{r_s}\log(y)\big).
\end{eqnarray*}
Using \eqref{beta1=1} we have
\[
\lim_{s\to\infty}F_{W_{n_s,\gamma_{n_s}}}(y)=\lim_{s\to\infty}F^*_{\underline{\boldsymbol\beta}}
(\frac{1}{r_s}\log(y))=F^*_{\underline{\boldsymbol\beta}}(0)=\frac12,
~~~y>0
\]
since $r_s\to\infty$. This means that
$\overline{W}_{n_s,\gamma_{n_s}}$ converges to infinity with
probability $0.5$, contradictory to the assumption that
$\overline{W}_{n_s,\gamma_{n_s}}$ converges weakly. Therefore, we
have proved that $\Delta_n/\gamma_n$ is bounded.

If there exists a subsequence of $\{n\}$ along which
$\Delta_n/\gamma_n$ converges to zero, then $\nu$ must be a
degenerate probability measure at $1$ from Part (\textbf{i}). Now we
can conclude that
\begin{equation}\label{boundedcoef}
 c_1\le \Delta_n/\gamma_n\le c_2~~\mbox{ for } n\ge 1
\end{equation}
for some constant $c_1,c_2\in (0,\infty)$.

Now we are ready to show \eqref{limits}.   Assume \eqref{limits} is
not true, that is,  for a $j_0$  the limit of
$\delta_{n,j_0}/\gamma_n$ doesn't exist.  Then there exist two
subsequences of $\{n\}$, say $\{n'\}$ and $\{n''\}$, such that
\begin{equation}\label{twobetaj0}
\lim_{n'\to\infty}\frac{\delta_{n',j_0}}{\gamma_{n'}}=\beta_{j_0}\ne
\hat{\beta}_{j_0}=\lim_{n''\to\infty}\frac{\delta_{n'',j_0}}{\gamma_{n''}}.
\end{equation}
From \eqref{boundedcoef} and \eqref{boundedratio}, we have
$\frac{|\delta_{n,j}|}{\gamma_n}=\frac{\Delta_n}{\gamma_n}\frac{|\delta_{n,j}|}{\Delta_n}\le
c_2$ for $n\ge 1$ and $j\ge 1$.  By using the diagonal argument, for
each of the two subsequences $\{n'\}$ and $\{n''\}$, there exists
its further subsequence along which $\frac{\delta_{n,j}}{\gamma_n}$
has a limit for all $j\ge 1$. Without introducing additional
notation, we can simply assume these limits exist along $\{n'\}$ and
$\{n''\}$, respectively, that is, we assume
\[
\lim_{n'\to\infty}\frac{\delta_{n',j}}{\gamma_{n'}}=\beta_j,~~
\lim_{n''\to\infty}\frac{\delta_{n'',j}}{\gamma_{n''}}=\hat\beta_j
\]
for $j\ge 1$.  Note that \eqref{boundedcoef} ensures $\beta_1>0$ and
$\hat\beta_1>0$.  Set $\boldsymbol{\beta}=(\beta_1, \beta_2,
\cdots)$ and $\hat{\boldsymbol{\beta}}=(\hat\beta_1, \hat\beta_2,
\cdots)$.  We apply Lemma~\ref{sublimits-haar} to the two
subsequences and conclude that $F^*_{\boldsymbol{\beta}}(y)$ and
$F^*_{\hat{\boldsymbol{\beta}}}(y)$ are the same for all $y$. This
implies $f_{\boldsymbol{\beta}}(0)=f_{\hat{\boldsymbol{\beta}}}(0)$,
$f_{\boldsymbol{\beta}}(1)=f_{\hat{\boldsymbol{\beta}}}(1)$, and
$f^{-1}_{\boldsymbol{\beta}}(y)=f^{-1}_{\hat{\boldsymbol{\beta}}}(y)$
in $(f_{\boldsymbol{\beta}}(0), f_{\boldsymbol{\beta}}(1))$, and
consequently,
$f_{\boldsymbol{\beta}}(x)=f_{\hat{\boldsymbol{\beta}}}(x)$ for
$x\in (0,\frac12)$. Therefore, the two analytic functions
$f_{\boldsymbol{\beta}}(z)$ and $f_{\hat{\boldsymbol{\beta}}}(z)$ on
$\{z\in\mathbb{C}:  |z-\frac12|<\frac12\}$ are the same, and their
$j_0$-th derivatives at $z=\frac12$ are also the same. This leads to
$\beta_{j_0}=\hat{\beta}_{j_0}$, contradictory to \eqref{twobetaj0}.
This completes the proof for Part (\textbf{ii}). \eop

\vspace{20pt}


\bibliographystyle{plain}
\bibliography{changjiangqi2024}

\end{document}